\documentclass[12pt]{article}
\usepackage{amssymb}
\usepackage{amsmath}
\usepackage{amscd}
\usepackage[all]{xy}
\usepackage[french]{babel}

\usepackage{color}
\title{Sur le groupe de Brauer transcendant}
\author{Jean-Louis Colliot-Th\'el\`ene et Alexei N. Skorobogatov}
\date{21 juillet 2011}

\def\sB{{\mathcal B}}
\def\sC{{\mathcal C}}

\def\beq{\begin{equation} \label}

\def\sA{{\mathcal A}}
\def\sB{{\mathcal B}}
\def\sC{{\mathcal C}}

\newtheorem{theo}{Th\'eor\`eme}[section]
\newtheorem{prop}[theo]{Proposition}
\newtheorem{lem}[theo]{Lemme}
\newtheorem{cor}[theo]{Corollaire}
\newtheorem{defi}[theo]{D\'efinition}

\newcommand{\bthe}{\begin{theo}}
\newcommand{\ble}{\begin{lem}}
\newcommand{\bpr}{\begin{prop}}
\newcommand{\bco}{\begin{cor}}
\newcommand{\bde}{\begin{defi}}
\newcommand{\ethe}{\end{theo}}
\newcommand{\ele}{\end{lem}}
\newcommand{\epr}{\end{prop}}
\newcommand{\eco}{\end{cor}}
\newcommand{\ede}{\end{defi}}

\def\lra{\longrightarrow}
\def\res{{\rm res}}
\def\hom{{\rm hom}}
\def\num{{\rm num}}
\def\H{{\rm H}}
\def\NS{{\rm NS}}
\def\CH{{\rm CH}}
\def\Num{{\rm Num}}
\def \Gal {{\rm{Gal}}}
\def \Spec {{\rm{Spec}}}
\def\et{{\rm{\acute et}}}
\def \Pic {{\rm {Pic}}}

\def \Ga {\Gamma}
\def\ov{\overline}
\def \Br {{\rm{Br}}}
\def \to {{\rightarrow}}
\def \Ker {{\rm Ker}}
\def \Coker {{\rm Coker}}
\def \Hom {{\rm Hom}}

\def \Im {{\rm Im}}
\def \cores {{\rm cores}}

\def \Q {{\mathbb Q}}
\def \N {{\mathbb N}}
\def \Z {{\mathbb Z}}

\def \C {{\mathbb C}}

\def \G {{\mathbb G}}
\def \P {{\mathbb P}}

\begin{document}

\maketitle

{\bf R\'esum\'e}

\noindent Soit $X$ une vari\'et\'e projective et lisse sur un corps $k$ de
carac\-t\'e\-ristique z\'ero.
Le groupe de Brauer de $X$ s'envoie dans les invariants,
sous le groupe de Galois absolu de $k$,   du groupe de Brauer de la m\^eme
vari\'et\'e consid\'er\'ee sur une cl\^oture alg\'ebrique de $k$.
Nous montrons que le quotient est fini. Sous des hypoth\`eses suppl\'ementaires,
par exemple sur   un corps de nombres, nous donnons des estimations sur
l'ordre de ce quotient. L'accouplement d'intersection entre les groupes de diviseurs et de 1-cycles modulo
\'equivalence num\'erique joue ici un r\^ole important.

\bigskip

{\bf Abstract}

\noindent For a smooth and projective variety $X$ over a field $k$
of characteristic zero we prove the finiteness of the 
cokernel of the natural map from the Brauer group of $X$ to the 
Galois-invariant subgroup of the Brauer group of the same
variety over an algebraic closure of $k$. Under further conditions
on $k$, e.g. over number fields, we give estimates for 
the order of this cokernel. We
emphasise the r\^ole played by the exponent of the discriminant
groups of the intersection pairing between the groups of divisors
and curves modulo numerical equivalence.

\medskip

MSC-class :   14F22 ; 14G99, 14C25, 16K50

\section*{Introduction}

Soit $X$ une vari\'et\'e projective, lisse et g\'eom\'etriquement int\`egre 
 sur
un corps $k$ de caract\'eristique z\'ero.  
Soit  $\ov k$ une cl\^oture alg\'ebrique de $k$.
Soient $\Gamma=\Gal(\ov k/k)$ et $\ov X= X \times_k\ov k$.
On a l'application naturelle de groupes de Brauer :
$$ \Br(X) \to \Br(\ov X).$$
Le noyau de cette application, not\'e  $\Br_1(X)$,
est appel\'e groupe de Brauer {\it alg\'ebrique}
de $X$. L'image de cette application est appel\'e
groupe de Brauer {\it transcendant} de $X$. C'est un sous-groupe
du groupe des invariants $\Br(\ov X)^\Gamma$. On a donc l'inclusion
$$\Br(X)/\Br_1(X) \subset \Br(\ov X)^\Gamma.$$
On voudrait calculer ces groupes, en particulier en vue
de l'\'etude de l'obstruction de Brauer--Manin.
La double question suivante a \'et\'e
soulev\'ee  
dans  \cite{SZ1} et \cite{CS}.

{\it  Si $k$ est un corps de type fini sur $\Q$, chacun des
groupes $\Br(X)/\Br_1(X)$ et $\Br(\ov X)^\Gamma$ est-il un groupe fini? }

 Dans cet article nous montrons que cette double question se r\'eduit
 \`a une seule question.

 \bigskip
 
 Pour $X$ une vari\'et\'e projective, lisse, et  g\'eom\'etriquement int\`egre
 sur  un corps $k$ de caract\'eristique z\'ero,  nous montrons que le conoyau de
 l'application naturelle
 $$\alpha:\Br(X)\to\Br(\ov X)^\Ga$$
est  un groupe fini. Sous des hypoth\`eses suppl\'ementaires sur le corps de base~$k$,
par exemple sur un corps de nombres, nous donnons des estimations
pour l'exposant et l'ordre de ce groupe fini.

Notre principal outil est un complexe naturel
(voir le paragraphe \ref{basicseq})
$$\Br(X)\buildrel{\alpha}\over{\lra} \Br(\ov X)^\Ga 
\buildrel{\beta}\over{\lra}\H^2(k,\Pic(\ov X))$$
qui pour $X$ avec un point rationnel ou pour $k$ un corps
de nombres est une suite exacte. Nous \'etudions l'image (finie)
de $\beta$ par deux m\'ethodes diff\'erentes, qui
m\`enent \`a des estimations similaires mais non identiques.

 La premi\`ere m\'ethode fait l'objet du paragraphe \ref{one}.
 L'id\'ee principale est d'utiliser la fonctorialit\'e du complexe
 ci-dessus par rapport aux morphismes de $k$-vari\'et\'es
 et d'utiliser la trivialit\'e du groupe de Brauer des courbes
 sur un corps alg\'ebriquement clos (th\'eor\`eme de Tsen).
 Cela montre que la restriction de l'image de $\beta$ 
 \`a toute courbe ferm\'ee dans $X$ est nulle. Ceci m\`ene
 aux th\'eor\`emes  \ref{te1} et \ref{te2}.

 La deuxi\`eme m\'ethode utilise une remarque g\'en\'erale sur les
 diff\'eren\-tielles dans la suite spectrale des foncteurs compos\'es.
 Soit $\Br^0(\ov X)$ le sous-groupe divisible maximal de $\Br(\ov X)$,
 et soit  $\NS(\ov X)$ le groupe de N\'eron--Severi. Nous montrons que
 l'application compos\'ee
 $$\Br^0(\ov X)^{\Ga}\hookrightarrow 
\Br(\ov X)^{\Ga} \buildrel{\beta}\over{\lra} 
\H^2(k,\Pic(\ov X)) \to \H^2(k,\NS(\ov X)/_{\rm tors})$$
peut se lire comme l'homomorphisme de connexion associ\'e
\`a une certaine 2-extension naturelle de  $\Ga$-modules
obtenue \`a partir de la suite de Kummer (voir le corollaire
\ref{c1}). D'apr\`es un  th\'eor\`eme de Lieberman rappel\'e
au paragraphe  \ref{1.1}, \'equivalence num\'erique et \'equivalence
homologique co\"{\i}ncident sur les cycles alg\'ebriques de dimension 1.
(Pour les surfaces, le th\'eor\`eme de Lieberman est un r\'esultat
classique de Matsusaka.) Nous utilisons ce fait pour montrer \`a la
proposition \ref{d2}  que l'image de l'application compos\'ee ci-dessus
est annul\'ee par l'exposant de chacun des deux groupes discriminants
d\'efinis par l'accouplement d'intersection entre les groupes de diviseurs
et les groupes de 1-cycles sur $\ov X$, modulo \'equivalence num\'erique.
Les bornes pour le conoyau de $\alpha$ obtenues par cette m\'ethode
sont donn\'ees aux th\'eor\`emes \ref{t} et \ref{t1}.

Au paragraphe \ref{surfaces} nous donnons des applications
aux surfaces K3  et aux produits de deux courbes. Pour toute
telle surface avec un $k$-point, l'\'enonc\'e est particuli\`erement simple :
le conoyau de $\alpha$ est annul\'e par l'exposant du groupe discriminant
d\'efini par la forme d'intersection  sur $\NS(\ov X)$, voir
les propositions \ref{K3} et  \ref{p}.

T.~Szamuely a demand\'e si le r\'esultat de finitude du th\'eor\`eme  \ref{te1}
vaut encore pour les vari\'et\'es lisses quasi-projectives. Au paragraphe
\ref{open} nous donnons une r\'eponse affirmative lorsque le corps
de base $k$ est un corps de type fini sur $\Q$.

Ce travail a \'et\'e commenc\'e lors de  la conf\'erence
``Arithmetic of surfaces" qui s'est tenue au Centre Lorentz
\`a Leiden en Octobre 2010. Nous en remercions les organisateurs.
Nous remercions  L.~Illusie, B.~Kahn,  J.~Riou pour  leur aide concernant le paragraphe \ref{1.1},
et  T.~Szamuely pour sa question.

\medskip

\section{Pr\'eliminaires}

Soient $k$ un corps de caract\'eristique z\'ero, 
$\ov k$ une cl\^oture alg\'ebrique de $k$ et 
$\Gamma={\rm Gal}(\ov k/k)$ le groupe de Galois absolu de $k$.

Soit  $X$ une vari\'et\'e projective, lisse et g\'eom\'etriquement int\`egre
sur $k$, de dimension $d$.
Soit  $\ov X= X \times_{k}\ov k$.  

Pour  un groupe ab\'elien $A$ et $n>0 $ un entier, on note $A[n] \subset A$ le sous-groupe
des \'el\'ements annul\'es par $n$. Pour $\ell$ un nombre premier on note
$A\{\ell\} \subset A$ le sous-groupe de torsion $\ell$-primaire.

\subsection{Cycles alg\'ebriques} \label{1.1}

Pour tout entier $i$ avec $0\leq i\leq d$, soit
$\CH^i(\ov X)$  le groupe de Chow des
cycles de codimension $i$ sur $\ov X$,  c'est-\`a-dire le groupe
des combinaisons lin\'eaires \`a coefficients entiers de sous-vari\'et\'es
ferm\'ees irr\'eductibles de codimension~$i$ modulo l'\'equivalence
rationnelle.  Comme $X$ est lisse, on a  $ \Pic(\ov X)=\CH^1(\ov X)$.
Soit $\NS(\ov X)$ le groupe de  N\'eron--Severi de $\ov X$. C'est
  le quotient de $\Pic(\ov X)$ par son sous-groupe divisible maximal $\Pic^0(\ov X)$,
  groupe des $\ov k$-points de la vari\'et\'e de Picard de $X$.

Puisque  $X$ est projective, l'intersection d\'efinit une forme bilin\'eaire
$\Ga$-\'equivariante  
\begin{equation}
\CH^i(\ov X)\times\CH^{d-i}(\ov X)\to\Z. \label{int_pairing}
\end{equation}
Soit $N^i=\Num^i(\ov X)$ le groupe des cycles de codimension $i$ sur $\ov X$
modulo \'equivalence num\'erique. C'est le quotient de $\CH^i(\ov X)$
par le noyau (\`a gauche) de l'accouplement (\ref{int_pairing}).
\'Ecrivons $N_i=N^{d-i}$. Nous obtenons une forme bilin\'eaire 
 $\Ga$-\'equivariante
 \begin{equation}
N^i\times N_i\to\Z \label{num_pairing}
\end{equation}
dont les noyaux \`a gauche et \`a droite sont triviaux.
Pour tout $i\geq 0$ le groupe ab\'elien
$N^i$ est libre de type fini. Ceci r\'esulte de l'existence
d'une cohomologie de Weil, \`a coefficients dans un corps
de caract\'eristique z\'ero et  munie d'applications classe de cycle
pour lesquelles le cup-produit en cohomologie est
compatible avec l'intersection des cycles  \cite[Thm. 3.5, p.~379]{Kl}.

Pour $i=1$, l'accouplement  (\ref{num_pairing}) donne naissance
\`a la suite exact de 
 $\Ga$-modules
\begin{equation}
0\to N^1\to \Hom(N_1,\Z)\to D\to 0,
\label{x4}
\end{equation}
qui d\'efinit le 
 $\Ga$-module fini $D$.
Ce groupe est l'un des deux groupes discriminants associ\'es 
\`a l'accouplement
 $N^1\times N_1\to\Z$.

Pour tout  $i \geq 0$ on dispose des applications classe de cycle
$$\CH^{i}(\ov X)\lra \H^{2i}_{\et}(\ov X,\Z_{\ell}(i)),$$
voir \cite[Section VI.9]{EC} et \cite[Cycle]{SGA41/2}.
Ces applications transforment
 cup-produit  en cohomologie $\ell$-adique en 
intersection des cycles alg\'ebriques, voir \cite[Prop. VI.9.5]{EC}.
Introduisons les $\Ga$-modules
$$N^i_\ell=N^i\otimes\Z_\ell,\quad N_{i,\ell}=N_i\otimes\Z_\ell,
\quad 
H^{2i}_\ell=\H_{\et}^{2i}(\ov X,\Z_\ell(i))/_{\rm tors}.$$
Pour toute th\'eorie cohomologique
 \`a coefficients
dans un corps de caract\'e\-ris\-tique z\'ero, 
avec applications classe de cycle
compatibles avec le cup-produit en cohomologie,
l'\'equivalence homologique implique l'\'equivalence
num\'erique. Selon les ``conjectures standards'', pour toute
bonne th\'eorie cohomologique, \'equivalence homologique
et \'equivalence num\'erique devraient co\"{\i}ncider.
Dans le contexte de la cohomologie $\ell$-adique,
l'application classe de cycle devrait se factoriser de la fa\c con suivante :
\begin{equation}\CH^{i}(\ov X) \otimes \Z_{\ell} \to N^i_\ell 
\hookrightarrow H^{2i}_{\ell}.\label{eq}
\end{equation}
C'est le cas pour $i=1$ par un th\'eor\`eme classique de
T. Matsusaka \cite{Mat}, qui montra $N^1=\NS(\ov X)/_{\rm tors}$.

Pour la cohomologie de Betti, avec applications classe de cycle
$$\CH^{i}(X) \to \H^{2i}_{\rm Betti}(X(\C),\Q(i)),$$
o\`u  $\Z(i)=\Z(2\pi \sqrt{-1})^{\otimes i}$, ceci fut
\'etabli 
pour $i=d-1$ par  D.~Lieberman \cite[Cor. 1]{Lie}. 
Une version plus alg\'ebrique de la d\'emonstration fut donn\'ee
par Kleiman \cite[Remark 3.10]{Kl}.  Ces deux articles 
\'etablissent des r\'esultats pour d'autres valeurs de $i$, et pour cela font appel
au th\'eor\`eme de l'indice de Hodge. Le cas
case $i=d-1$ est plus simple, comme nous expliquons maintenant.

\begin{prop} \label{num1=hom1}  
Soit $X$ une vari\'et\'e connexe, projective et lisse sur $\C$.
 D\'efi\-nissons l'\'equivalence homologique sur
les cycles au moyen de la cohomologie de Betti \`a coefficients
rationnels. L'homomorphisme naturel 
$$CH_{1}(X)/\hom \lra CH_{1}(X)/\num$$
est un isomorphisme.
\end{prop}
{\it D\'emonstration.}
On peut supposer $d=\dim(X) \geq 3$.
Soit $L \in CH^1(X)$ la classe d'une section hyperplane.
Pour  $A$ un groupe ab\'elien, on note  $A_{\Q}:=A \otimes_{\Z}\Q$.
La multiplication par
 $L^{d-2} \in CH^{d-2}(X)$ d\'efinit un diagramme commutatif d'espaces
 vectoriels sur $\Q$ :
$$\begin{array}{cccccccccccccccccc}
CH^{d-1}(X)_{\Q}/\num  & \leftarrow & CH^{d-1}(X)_{\Q}/\hom  & \hookrightarrow &  Hdg^{d-2}(X,\Q) &   \hookrightarrow  & H^{2d-2}(X,\Q(d-1))\\
\uparrow  & &\uparrow  &&\uparrow  && \uparrow  \\
CH^{1}(X)_{\Q}/\num  & \leftarrow & CH^{1}(X)_{\Q}/\hom  &  \hookrightarrow &  Hdg^{2}(X,\Q) &  \hookrightarrow & H^{2}(X,\Q(1)) .\\
\end{array}$$

Les fl\`eches horizontales pointant vers la gauche sont surjectives.
Toutes les fl\`eches horizontales pointant vers la droite sont par d\'efinition injectives.
 D'apr\`es le th\'eor\`eme de Lefschetz difficile,  la quatri\`eme fl\`eche verticale est un 
isomorphisme. La d\'ecomposition de Hodge des groupes 
$H^{i}(X,\C)$ et le fait que les classes de type $(p,q)$ s'envoie sur des classes de type $(p+d-2,q+d-2)$ 
pour $p=0,1,2$ implique alors que la troisi\`eme fl\`eche verticale, qui porte sur les classes de Hodge, est  
un isomorphisme.
Par le th\'eor\`eme de Lefschetz sur les classes de type $(1,1)$, 
l'application $CH^{1}(X)_{\Q}/\hom   \rightarrow    Hdg^{2}(X,\Q)$ est un isomorphisme.
 Tout ceci implique que la deuxi\`eme fl\`eche verticale est aussi un isomorphisme.
 La fl\`eche verticale de gauche est donc surjective.
 Par d\'efinition, les deux espaces vectoriels de dimension finie
 $CH^{1}(X)_{\Q}/\num$ et $CH^{d-1}(X)_{\Q}/\num$
 ont la m\^{e}me dimension. La fl\`eche verticale de gauche est donc un isomorphisme.
D'apr\`es le th\'eor\`eme de Matsusaka,  la fl\`eche inf\'erieure gauche est un 
isomorphisme. On conclut que l'application $$CH^{d-1}(X)_{\Q}/\hom \to CH^{d-1}(X)_{\Q}/\num$$
est un isomorphisme. Ceci implique que l'application $$CH^{d-1}(X)/\hom \to CH^{d-1}(X)/\num$$
est un isomorphisme de groupes ab\'eliens de type fini sans torsion.
 QED

\bigskip

Rappelons maintenant comment divers th\'eor\`emes de comparaison
impliquent (\ref{eq})
pour $i=d-1$, 
o\`u $X$ est une vari\'et\'e projective, lisse, g\'eom\'e\-tri\-quement int\`egre sur
un corps $k$ de caract\'eristique z\'ero. Rappelons que pour un corps
alg\'ebriquement clos $L$ contenant $k$ le groupe de N\'eron--Severi
de $X_L=X\times_kL$  ne d\'epend pas du corps $L$, car c'est le groupe
des composantes connexes du sch\'ema de Picard  $\Pic_{X_L/L}$.
Nous pouvons donc utiliser la notation $N^1$  sans risque d'ambigu\"{\i}t\'e.

 Soit $C$ un 1-cycle sur $\ov X$ qui est num\'eriquement \'equivalent
 \`a z\'ero. Il existe un sous-corps $K\subset \ov k$ de type fini sur $\Q$,
 une vari\'et\'e $\tilde X$ sur $K$, et un  1-cycle $\tilde C$ sur  $\tilde X$
tel que $X=\tilde X\times_K\ov k$ et $C=\tilde C\times_K\ov k$. 
On peut supposer que le groupe de type fini $N^1$ est engendr\'e par les classes
de diviseurs effectifs, r\'eduits, absolument irr\'eductibles
$D_1,\ldots,D_r$ d\'efinis sur $K$. Choisissons un plongement
$K\subset \C$. Soit $\ov K$ la cl\^{o}ture alg\'ebrique de  $K$ dans $\C$,
et soient $\tilde X_{\ov K}=\tilde X\times_K\ov K$,
$\tilde X_{\C}=\tilde X\times_K\C$. Le cycle
 $C$ a une image nulle dans $N_1=\Num_1(\ov X)$
si et seulement si  $C$ a une intersection nulle avec  $D_1,\ldots,D_r$. 
Mais alors  $\tilde C$ a une image nulle dans  $\Num_1(\tilde X_\C)$. 
D'apr\`es la  proposition  \ref{num1=hom1},  le cycle  
$\tilde C$ 
a une image nulle dans le groupe de cohomologie de Betti  $\H^{2d-2}(\tilde X_\C,\Q(d-1))$.

Le th\'eor\`eme de comparaison entre la cohomologie  \'etale et la cohomologie de Betti
(\cite[XI, XVI]{SGA4},   voir aussi \cite[Thm. III.3.12]{EC})
donne des isomorphismes naturels
$$\H^{2i}_\et(\tilde X_\C,\Q_\ell(i))\cong \H^{2i}_{\rm Betti}(\tilde X_\C(\C),\Q(i))\otimes_\Q\Q_\ell.$$
Les applications classe de cycle transforment cup-produit en cohomologie en
accouplement d'intersection sur les groupes de Chow.

Partant de cela, on peut montrer que l'application classe de cycle  en
cohomologie de Betti et l'application classe de cycle en cohomologie
$\ell$-adique sont compatibles avec ces isomorphismes.
Une esquisse de d\'emonstration est donn\'ee dans
\cite[p.~21]{LNM900}. J.~Riou  nous a montr\'e comment
une preuve formelle se d\'eduit de l'\'enonc\'e d'unicit\'e pour les
applications classe de cycle que l'on trouve dans
  \cite[Prop. 1.2]{Riou}. 
 

Puisque l'application naturelle
$$\H^{2d-2}_\et(\tilde X_{\ov K},\Q_\ell(d-1))\lra\H^{2d-2}_\et(\tilde X_\C,\Q_\ell(d-1))$$
est un isomorphisme d'espaces vectoriels sur $\Q_\ell$ (cf. \cite[Cor. VI.4.3]{EC}),
l'application  classe de cycle 
envoie  $\tilde C$  sur z\'ero dans
$\H^{2d-2}_\et(\tilde X_{\ov K},\Q_\ell(d-1))$. 
Par changement de corps de base
de $\ov K$ \`a $\ov k$, on obtient (\ref{eq}) for $i=d-1$.

\medskip

Comme rappel\'e ci-dessus, l'accouplement
 (\ref{int_pairing})
est compatible avec le cup-produit
$$\H_{\et}^{2i}(\ov X,\Z_\ell(i)) \times
\H_{\et}^{2d-2i}(\ov X,\Z_\ell(d-i))\to\Z_\ell$$
via l'application 
classe de cycle.
Nous obtenons donc le diagramme commutatif d'accouplements
de
$\Ga$-modules
\begin{equation}
\begin{array}{ccccc}
N^1& \times& N_1 & \to& \Z\\
\downarrow&&\downarrow&&\downarrow\\
H^2_\ell&\times&H^{2d-2}_\ell&\to&\Z_\ell
\end{array}\label{pd}\end{equation}
o\`u les applications verticales sont injectives.
Nous utiliserons l'\'enonc\'e suivant :
l'accouplement inf\'erieur dans  (\ref{pd}) 
est un accouplement parfait, c'est-\`a-dire
qu'il induit des isomorphismes
$$H^2_\ell=\Hom_{\Z_\ell}(H^{2d-2}_\ell,\Z_\ell), \quad
H^{2d-2}_\ell=\Hom_{\Z_\ell}(H^2_\ell,\Z_\ell).$$
L.~Illusie nous informe que cet  \'enonc\'e
peut \^etre \'etablie en utilisant le formalisme $\Z_{\ell}$-adique
de Deligne \cite[\S 1.1]{DeligneW2}. La dualit\'e de Poincar\'e
pour le complexe $R{\Gamma}(X,\Z/\ell^n)$ (voir \cite[XVIII]{SGA4})
donne naissance \`a une dualit\'e parfaite pour les  complexes parfaits
$R{\Gamma}(X,\Z_{\ell})$. On utilise ensuite un argument 
de type coefficients universels.

La suite (\ref{x4}) donne naissance \`a la suite exacte
\begin{equation}
0\to N^1_\ell\to \Hom_{\Z_\ell}(N_{1,\ell},\Z_\ell)\to D\{\ell\}\to 0,
\label{x6}
\end{equation}
o\`u la seconde fl\`eche se factorise de la fa\c con suivante :
$$N^1_\ell\to H^2_\ell\tilde\lra\Hom_{\Z_\ell}(H^{2d-2}_\ell,\Z_\ell)
\to\Hom_{\Z_\ell}(N_{1,\ell},\Z_\ell).$$

\subsection{Le groupe de Brauer}\label{1.2}

Rappelons le calcul du groupe de Brauer $\Br(\ov X)$
(Grothendieck,  \cite[III.8, p.~144-147]{G68}).
Soit  $\rho={\rm dim}_\Q (\NS(\ov X)\otimes \Q)$  le nombre de Picard de
$\ov X$,  et soit
$b_2$ le second nombre de Betti de $\ov X$.
Notons $\Br^0(\ov X)$ le sous-groupe divisible maximal de
 $\Br(\ov X)$.
 On a un isomorphisme de groupes ab\'eliens :
 $$\Br^0(\ov X)\cong(\Q/\Z)^{b_2-\rho}.$$
Le quotient $\Br(\ov X)/\Br^0(\ov X)$ est fini, 
plus pr\'ecis\'ement il y a une suite exacte de
 $\Ga$-modules
\begin{equation}
0 \to \Br^0(\ov X) \to \Br(\ov X) \to 
\oplus_{\ell}\H^3_{\et}(\ov X,\Z_{\ell}(1))_{\rm tors} \to 0,
\label{brgeneral}
\end{equation}
o\`u $\ell$ parcourt l'ensemble des nombres premiers.

Soit $B_\ell$ le module de Tate  $\ell$-adique de $\Br(\ov X)$, que l'on d\'efinit
comme la limite projective des $\Br(\ov X)[\ell^m]$, $m \in \N$.
 C'est un  $\Z_\ell$-module libre de type fini.
Le module galoisien $B_\ell$ ne contr\^{o}le que le sous-groupe divisible
maximal
$\Br^0(\ov X)\subset \Br(\ov X)$, en ce sens que $B_\ell$
est aussi isomorphe au module de Tate de $\Br^0(\ov X)$, 
et qu'il y a un isomorphisme canonique de
$\Ga$-modules (cf. 
\cite[II.8.1, p.~144]{G68}) :
$$\Br^0(\ov X)\cong \oplus_\ell(B_\ell\otimes_{\Z_\ell}\Q_\ell/\Z_\ell).$$

La suite de Kummer
\begin{equation}
1\to\mu_{n}\to\G_m \buildrel{x \mapsto x^n}\over{\longrightarrow} \G_m\to 1 \label{kum}
\end{equation}
donne naissance aux suites exactes de $\Ga$-modules
$$0\to \Pic(\ov X)/\ell^m\to \H_{\et}^2(\ov X,\mu_{\ell^m})
\to \Br(\ov X)[{\ell^m}] \to 0.$$
Puisque le groupe divisible  $\Pic^0(\ov X)$ a une image nulle dans
$\H_{\et}^2(\ov X,\mu_{\ell^m})$, on a des suites exactes induites :
$$0\to \NS(\ov X)/\ell^m\to \H_{\et}^2(\ov X,\mu_{\ell^m})
\to \Br(\ov X)[{\ell^m}]\to 0.$$
En passant \`a la limite projective sur  $m \in \N$, on obtient la
suite exacte
  (8.7) de \cite[III.8.2]{G68} :
\begin{equation}
0\to \NS(\ov X)\otimes \Z_{\ell }\to \H_{\et}^2(\ov X,\Z_\ell(1))
\to B_\ell\to 0.\label{x1}
\end{equation}
La deuxi\`eme fl\`eche dans (\ref{x1})
induit un isomorphisme sur les groupes de torsion :
$$(\NS(\ov X)\otimes \Z_{\ell})_{\rm tors}=
\H_{\et}^2(\ov X,\Z_\ell(1))_{\rm tors}.$$
On a donc la suite exacte de $\Z_{\ell}[\Ga]$-modules,
libres et de type fini comme $\Z_{\ell}$-modules :
\begin{equation}
0\to N^1_\ell\to H^2_\ell\to B_\ell\to 0. \label{tate}
\end{equation}
Comme suite de $\Z_{\ell}$-modules,  cette suite est scind\'ee.
En particulier, pour tout premier~$\ell$ le $\Z_\ell$-sous-module 
$N^1_\ell\subset H^2_\ell$
est primitif, en ce sens que le quotient $ H^2_\ell/N^1_{\ell}$ est sans torsion.

Tensorisant (\ref{tate}) avec $\Q_\ell/\Z_\ell$ et prenant la somme directe sur tous
les premiers  $\ell$, nous obtenons une suite exacte de $\Ga$-modules
$$0\to N^1\otimes\Q/\Z\to \oplus_\ell (H^2_\ell\otimes_{\Z_\ell}\Q_\ell/\Z_\ell)
\to \Br^0(\ov X)\to 0,$$
qui donne naissance \`a une  2-extension de $\Ga$-modules
\begin{equation}
0\to N^1\to N^1\otimes\Q\to \oplus_\ell (H^2_\ell\otimes_{\Z_\ell}\Q_\ell/\Z_\ell)
\to \Br^0(\ov X)\to 0. \label{miau}
\end{equation}

\medskip

On utilisera plus loin le lemme facile suivant.

\ble\label{trivial}
Soit $F$ un groupe ab\'elien  fini $\ell$-primaire, et soit $n \in \N$.
Soit  $A$ un sous-quotient fini de $(\Q_{\ell}/\Z_{\ell})^n \oplus F$.
Si l'exposant de  $A$ est  $\ell^m$, alors l'ordre de $A$
divise le produit de $\ell^{mn}$ par l'ordre de $F[\ell^m]$.
\ele
{\it D\'emonstration.} Le groupe  $A$ est un quotient de
$(\Q_{\ell}/\Z_{\ell})^r \oplus F' \subset 
(\Q_{\ell}/\Z_{\ell})^n \oplus F$,
o\`u $F'$ est un groupe fini. Donc  $A$ est un quotient de $F'/\ell^m$.
L'ordre de  $F'/\ell^m$ est \'egal \`a l'ordre de $F'[\ell^m]$, 
qui est un sous-groupe de
$(\Z/\ell^m)^n \oplus F[\ell^m]$. QED

\subsection{Une suite exacte fondamentale}\label{basicseq}

\bpr \label{bas}
Soit $X$ un sch\'ema de type fini sur un corps $k$
de caract\'eristique z\'ero.

{\rm (i)} Il y a un complexe naturel, fonctoriel en $X$ et  en $k$:
$$\Br(X)\buildrel{\alpha}\over{\lra} \Br(\ov X)^\Ga 
\buildrel{\beta}\over{\lra}\H^2(k,\Pic(\ov X)).$$

{\rm (ii)} Supposons $\H^0_\et(\ov X,\G_m)=\ov k^*$.
Supposons de plus que l'application
$\H^3_\et(k,\ov k^*) \to  \H^3_{\et}(X,\G_{m})$est injective,
ce qui est le cas si $X$ poss\`ede un $k$-point ou si  $k$ est un corps
de nombres.
Alors le complexe ci-dessus est une suite exacte, et l'on a
 $\Im(\alpha)=\Ker(\beta)$ et $\Coker(\alpha)= \Im(\beta)$. 
\epr
{\it D\'emonstration.}
Ceci r\'esulte de la suite spectrale de Leray
\begin{equation}
E_{2}^{pq}=\H^p(k,\H^q_{\et}(\ov X,\G_m))\Rightarrow \H^{p+q}_{\et}(X,\G_m).
\label{ss1}
\end{equation}
Un  $k$-point sur  $X$ d\'efinit une section de l'application
$\H^3_\et(k,\ov k^*) \to  \H^3_{\et}(X,\G_{m})$.
Pour un corps de nombres $k$, on a $\H^3_\et(k,\ov k^*)=0$. QED
 
\subsection{Restriction et corestriction}

Le lemme suivant est certainement bien connu.

\ble \label{comm}
Soit $X$ un sch\'ema sur un corps $k$ de caract\'eristique nulle,
et soit
 $L\subset \ov k$ une extension finie de  $k$ de degr\'e  $n$. 
Il existe des homomorphismes naturels de restriction et de corestriction
$$\res_{L/k}:\Br(X)\to\Br(X_L), \quad\quad\cores_{L/k}:\Br(X_L)\to\Br(X),$$
et  l'on  a $\cores_{L/k}(\res_{L/k}(x))=nx$. Le diagramme suivant commute :
$$\xymatrix{
\Br(X)\  \ar[d]_\alpha \ar[r]^{\res_{L/k}}& \
\Br(X_L)\ \ar[d]_{\alpha_L} \ar[r]^-{\cores_{L/k}}&
\ \Br(X) \ar[d]_\alpha\\
\Br(\ov X)^\Ga\ \ar@{^{(}->}[r] &\ \Br(\ov X)^{\Ga_L}\ 
\ar[r]^\sigma &\ \Br(\ov X)}
$$
Ici $\Ga_L=\Gal(\ov k/L)$, et $\sigma(x)=\sum \sigma_i(x)$, o\`u
les $\sigma_i\in \Ga$ sont des repr\'esentants des classes $\Ga/\Ga_L$. 
\ele


{\em D\'emonstration.} 
Rappelons la d\'efinition de $\res_{L/k}$  et $\cores_{L/k}$. 
Soit $f:Y\to X$ un morphisme fini et plat de $k$-sch\'emas lisses
connexes. Soit $n$ le degr\'e de $f$. On dispose alors de morphismes
de faisceaux \'etales $$\G_{m,X}\to f_*\G_{m,Y}\to \G_{m,X}$$ d\'efinis
sur les fibres par l'injection naturelle pour le premier, par la norme pour
le second. Le morphisme compos\'e est l'\'el\'evation \`a la puissance $n$.
Le foncteur $f_*$  de la cat\'egorie des faisceaux \'etales sur $Y$ dans
la cat\'egorie des faisceaux \'etales sur $X$ est exact \cite[Cor. II.3.6]{EC}.
La suite spectrale de Leray donne donc un isomorphisme
$\H^p_\et(X,f_*\G_{m,Y})\tilde\lra \H^p_\et(Y,\G_{m,Y})$. 
On obtient ainsi les applications 
$$
\H^p_\et(X,\G_{m,X})\buildrel{\res}\over{\hbox to 16mm{\rightarrowfill}}
\H^p_\et(Y,\G_{m,X})\buildrel{\cores}\over{\hbox to 16mm{\rightarrowfill}}
\H^p_\et(X,\G_{m,X})
$$
dont la compos\'ee est la multiplication par $n$.

Soit $X$ un sch\'ema sur un corps $k$. Soit $L\subset \ov k$
un corps tel que $[L:k]=n$. 
Soit  $Y=X_L=X\times_kL$.
On a l'isomorphisme $L\otimes_k\ov k  \tilde\to \ov k^n$, 
dont les diverses composantes correspondent aux $n$ $k$-plongements
de $L$ dans $\ov k$.

Par changement de base de $X$ \`a $\ov X$, on obtient  un diagramme commutatif
$$\xymatrix{
\H^p_\et(X,\G_m)\  \ar[d] \ar[r]^{\res_{L/k}}& \
\H^p_\et(X_L,\G_m)\ \ar[d] \ar[r]^-{\cores_{L/k}}&
\ \H^p_\et(X,\G_m) \ar[d]\\
\H^p_\et(\ov X,\G_m)\ \ar@{^{(}->}[r] &\ \H^p_\et(\ov X,\G_m)^n\ 
\ar[r]&\ \H^p_\et(\ov X,\G_m)}
$$
o\`u les applications dans la ligne inf\'erieure sont le plongement diagonal
et le produit. L'action du groupe de Galois $\Ga$ sur $\H^p_\et(\ov X,\G_m)^n$
est induite par l'action naturelle  de $\Ga_L$  sur $\H^p_\et(\ov X,\G_m)$.
En passant aux sous-groupes $\Ga$-invariants, et en prenant $p=2$,
on obtient l'\'enonc\'e du lemme. QED



\section{D\'emonstration du  th\'eor\`eme principal via le groupe de Brauer des courbes} \label{one}

\subsection{Finitude}

\bthe \label{te1}
Soit $X$ une vari\'et\'e projective, lisse et g\'eom\'etriquement int\`egre sur
un corps $k$ de caract\'eristique z\'ero. Le conoyau de l'application naturelle
 $\alpha:\Br(X)\to \Br(\ov X)^\Gamma$
est fini.
\ethe
{\it D\'emonstration.} D'apr\`es le calcul de Grothendieck du groupe
$\Br(\ov X)$ que nous avons rappel\'e au paragraphe  \ref{1.2}, 
pour toute puissance $\ell^n$ d'un nombre premier $\ell$
et tout sous-quotient $B$ de $\Br(\ov X)$ le sous-groupe $B[\ell^n]$ 
est fini. Il suffit donc de montrer que le groupe
$\Coker(\alpha)$  est d'exposant fini.

Pour cela on peut remplacer $k$ par une extension finie.
De fait, si $L$ avec $k\subset L\subset\ov k$, $[L:k]=n$, 
est une telle extension, il r\'esulte du lemme  \ref{comm}
que l'on a des applications naturelles 
 $$\Coker(\alpha)\to\Coker(\alpha_L) \to\Coker(\alpha),$$
 dont la compos\'ee est la multiplication par $n$.
 Il suffit donc de montrer que $\Coker(\alpha_L)$
est d'exposant fini.

On peut en particulier supposer que $X$ poss\`ede un $k$-point.
D'apr\`es la proposition 
  \ref{bas}(ii)
on a  $\Coker(\alpha)=\Im(\beta)$.
Montrons que $\Im(\beta)$ est d'exposant fini.

Si $C$ est une courbe projective, lisse et  g\'eom\'etriquement int\`egre
sur $k$ et si $f:C\to X$  est un $k$-morphisme, les applications
$f^*:\Pic(\ov X)\to \Pic(\ov C)$  et $f^*:\Br(\ov X)\to \Br(\ov C)$
s'ins\`erent dans  un diagramme commutatif  
$$\begin{array}{ccc}
\Br(\ov X)^\Ga&\buildrel{\beta_{X}}\over{\lra} &\H^2(k,\Pic(\ov X))\\
\downarrow&&\downarrow\\
\Br(\ov C)^\Ga&\buildrel{\beta_{C}}\over{\lra} &\H^2(k,\Pic(\ov C)).
\end{array}$$
Le th\'eor\`eme de Tsen implique $\Br(\ov C)=0$
(\cite{G68}, Cor. 1.3, p. 90).
Ainsi

\smallskip

  $(*)$
{\it  Pour tout $k$-morphisme  $f : C \to X$, le groupe $\Im(\beta_{X})$
est dans le noyau de l'application verticale droite du diagramme ci-dessus. }
 
\smallskip

L'application degr\'e $\Pic(\ov C)\to\NS(\ov C)=\Z$ donne naissance \`a la suite exacte
de modules galoisiens
$$0\to \Pic^0(\ov C)\to \Pic(\ov C)\to\NS(\ov C)\to 0.$$
On a donc un diagramme commutatif \`a lignes exactes
\begin{equation}
\begin{array}{ccccccc}
&&\H^2(k,\Pic^0(\ov X))&\to &\H^2(k,\Pic(\ov X))&\to
&\H^2(k,\NS(\ov X))\\
&&\downarrow&&\downarrow&&\downarrow \\
0&\to&\H^2(k,\Pic^0(\ov C))&\to &\H^2(k,\Pic(\ov C))&\to&\H^2(k,\NS(\ov C)).
\end{array}\label{di1}
\end{equation}
La nullit\'e de $\H^1(k,\Z)$ donne le z\'ero \`a gauche dans la ligne inf\'erieure.

Comme le corps $k$ est infini, le th\'eor\`eme de Bertini \cite{Jou}
pour les sections hyperplanes des vari\'et\'es projectives et lisses
implique l'existence d'une courbe $C \subset X$, definie sur  $k$,
section lin\'eaire de $X$, et qui est lisse et g\'eom\'etriquement connexe.
Une combinaison du th\'eor\`eme de Bertini  et du th\'eor\`eme de connexion
de Zariski (voir \cite[Lemme~2.10, p.~210]{SGA1})  
montre alors que sur une cl\^{o}ture alg\'ebrique de $k$,
l'image inverse via $f : C \to X$ de tout rev\^etement fini
\'etale connexe de $X$ est connexe. Ceci implique en
particulier que l'homomorphisme de vari\'et\'es ab\'eliennes
$\Pic^0_{X/k} \to \Pic^0_{C/k}$ a un noyau trivial.

 Par le th\'eor\`eme de compl\`ete r\'eductibilit\'e de Poincar\'e
 \cite[\S 19, Thm.~1]{Mu}
il existe donc une sous-var\'et\'e ab\'elienne  $A \subset \Pic^0_{C/k}$
telle que l'application naturelle
$$\Pic^0_{X/k} \times A \to \Pic^0_{C/k}$$
soit une isog\'enie de vari\'et\'es ab\'eliennes sur $k$.

Puisque  l'application  $\H^2(k,\Pic^0(\ov C)) \to \H^2(k,\Pic(\ov C))$ est  injective,
ceci implique :
 
\smallskip

$(**)$ {\it Le noyau de l'application compos\'ee
$$\H^2(k, \Pic^0({\ov X})) \to \H^2( k, \Pic({\ov X})) \to 
\H^2(k, \Pic({\ov C}))$$
a un exposant fini.}

\smallskip

Puisque  $N^1=\NS(\ov X)/_{\rm tors}$ est un groupe ab\'elien libre de type fini,
on peut choisir un nombre fini, disons $m$, de courbes int\`egres  sur $\ov X$ 
telles que l'intersection avec ces courbes d\'efinit un homomorphisme injectif
$\iota : N^1 \hookrightarrow \Z^m$. Par passage aux normalisations
on obtient des  morphismes  
de courbes projectives, lisses, connexes
d\'efinies sur $\ov k$, vers $\ov X$.  Pour la pr\'esente d\'emonstration
on peut remplacer $k$ par une extension finie sur laquelle chacune
des courbes  est d\'efinie.

Nous avons donc des $k$-courbes $C_{i}$, $i=1,\ldots,m$, projectives, lisses, g\'eo\-m\'e\-triquement
int\`egres et des $k$-morphismes $f_i:C_i\to X$. Les applications induisent un
homomorphisme de $\Gamma$-modules
$$   \NS(\ov X) \to \oplus_{i=1}^m \NS({ \ov C}_{i})= \Z^m.$$
D'apr\`es (\ref{di1}), $(*)$ et $(**)$, 
pour \'etablir le r\'esultat annonc\'e, il suffit de montrer que le
noyau de l'application induite
$\H^2(k,\NS(\ov X))\to \H^2(k,\Z^m)$ 
est d'exposant fini.
Cette application est la compos\'ee de deux applications :
$$\H^2(k,\NS(\ov X))\to \H^2(k,N^1)  \to \H^2(k,\Z^m).$$
Il suffit   de montrer  que le noyau de chacune de ces deux applications
est d'exposant fini.

 De la suite exacte de cohomologie associ\'ee \`a la suite exacte de 
 $\Ga$-modules
$$0\to \NS(\ov X)_{\rm tors}\to \NS(\ov X)\to N^1 \to 0$$
on d\'eduit que l'application
$$\H^2(k, \NS(\ov X))\to \H^2(k,N^1)$$
a son noyau annul\'e par la multiplication par l'exposant du groupe fini
$\NS(\ov X)_{\rm tors}$.
Il existe un homomorphisme $\Z^m  \to  N^1$ 
tel que la composition d'homomorphismes de groupes ab\'eliens,
avec action triviale du groupe de Galois,
$$N^1  \buildrel{\iota}\over{\hookrightarrow}  \Z^m \to N^1$$
est la multiplication par un entier strictement positif.
Le noyau de l'application  $$\H^2(k,N^1) \to \H^2(k,\Z^m)$$
est annul\'e par la multiplication par cet entier. QED

\subsection{Majorations, I}

Soit  $\delta_0$ l'exposant du groupe fini $D$ 
d\'efini dans (\ref{x4}), et soit
$\nu_0$ l'expo\-sant du groupe fini  $\NS(\ov X)_{\rm tors}$. 
Soit $\alpha$ l'application naturelle $\Br(X)\to \Br(\ov X)^\Gamma$.

\bthe \label{te2}
Soit  $X$ une vari\'et\'e projective, lisse et  g\'eom\'etriquement in\-t\`e\-gre
sur un corps $k$ de caract\'eristique z\'ero.
Soit $L/k$ une extension finie telle que le groupe ab\'elien 
libre de type fini
$N_1=\Num_1(\ov X)$
est engendr\'e par les classes de courbes int\`egres sur 
 $\ov X$ d\'efinies sur
 $L$. Let $\lambda = [L:k]$.

{\rm (i)} Si l'on a $\H^1(X,O_{X})=0$ et si l'application
$\H^3_{\et}(k,\G_{m}) \to \H^3_{\et}(X,\G_{m})$ 
est injective, alors l'exposant de  $\Coker(\alpha)$ 
divise  $\lambda \delta_0\nu_0$.

{\rm (ii)} Si $k$ est un corps de nombres, l'exposant de   $\Coker(\alpha)$
divise $2\lambda \delta_0\nu_0$, et il divise $\lambda \delta_0\nu_0$ 
si $k$ est totalement imaginaire.
\ethe
{\it D\'emonstration.} Elle consiste \`a d\'etailler les \'etapes
de la d\'emonstration du th\'eor\`eme \ref{te1}.
On applique d'abord la proposition  \ref{bas}.
Pour tout corps de nombres $k$, on a  $\H^3_{\et}(k,\ov k^*)=0$.
Dans ce cas,  on a donc  
 $\Im(\alpha)=\Ker(\beta)$,
et donc  $\Coker(\alpha)=\Im(\beta)$.
Ceci vaut aussi sous l'hypoth\`ese que l'application
$\H^3_{\et}(k,\G_{m}) \to \H^3_{\et}(X,\G_{m})$ est injective.

On choisit un nombre fini, soit  $m$, de courbes int\`egres
  $C_{1}, \dots, C_{m}$ sur $\ov X$ 
dont les classes engendrent $N_1$, et on remplace
 $k$ par  une extension finie $L$
 sur laquelle chacune de ces courbes est d\'efinie.
 L'argument de 
restriction-corestriction au d\'ebut de la d\'emonstration du th\'eor\`eme
 \ref{te1} montre que si l'on remplace
  $k$ par l'extension finie  $L$, de degr\'e  $\lambda$,
l'exposant $\Coker(\alpha)$ divise le produit de l'exposant de
  $\Coker(\alpha_{L})$ par l'entier  $\lambda$.
Pour \'etablir le th\'eor\`eme, il suffit donc de se limiter au cas
$k=L$,  c'est-\`a-dire  \`a $\lambda=1$.

Soit  $\Pic(\ov X) \to \Z^m$ l'application donn\'ee par restriction aux courbes 
  $C_{i}$, suivie de l'application degr\'e sur chaque courbe.
La d\'emonstration du th\'eor\`eme  \ref{te1} \'etablit  que l'image de  $\beta$
est contenue dans le noyau de l'application induite
$$\H^2(k,\Pic(\ov X)) \to \H^2(k,\Z^m).$$
L'application  $\Pic(\ov X) \to \Z^m$ se factorise comme suit :
$$\Pic(\ov X)\to\NS(\ov X)\to N^1\to 
\Hom(N_1,\Z) \to \Z^m.$$
Nous allons borner l'exposant 
du noyau 
de chaque application induite sur
 $\H^2(k,\bullet)$.

En envoyant chaque courbe  $C_{i}$ sur sa classe dans $N_1$
on obtient une suite exacte de $\Ga$-modules triviaux
$$0\to \Z^r\to \Z^m\to N_1\to 0.$$ 
En dualisant cette suite on obtient une suite exacte scind\'ee de
 $\Ga$-modules triviaux :
$$0\to \Hom(N_1,\Z) \to \Z^m\to\Z^r\to 0.$$
L'application
$$\H^2(k,\Hom(N_1,\Z)) \to \H^2(k, \Z^m)$$
est donc injective.
De la suite exacte  (\ref{x4}) on tire que le noyau de
  $$\H^2(k,N^1)\to \H^2(k,\Hom(N_1,\Z))$$
  est annul\'e par l'exposant de $\H^1(k,D)$, donc par $\delta_0$, 
  l'exposant du groupe $D$.
 
  Comme on a vu dans la d\'emonstration du th\'eor\`eme pr\'ec\'edent,
  le noyau de l'application
$$\H^2(k, \NS(\ov X))  \to \H^2(k,N^1)$$
est annul\'e par l'exposant
 $\nu_{0}$ de $\NS(\ov X)_{\rm tors}$.
Nous avons par ailleurs la suite exacte
$$ 0 \to \Pic^0_{X/k}(\ov k) \to \Pic(\ov X) \to \NS(\ov X) \to 0.$$
Si $\H^1(X,O_{X})=0$, alors $\Pic^0_{X/k}=0$.
Si  $k$ est un corps de nombres et  $A$ une vari\'et\'e ab\'elienne,
on  a  $\H^2(k,A)=\oplus_{v}\H^2(k_{v},A)$, o\`u  $v$ parcourt les places
r\'eelles de  $k$
\cite[Thm. 6.26 (c), p.~92]{Milne}. Ainsi l'exposant de
$\H^2(k,A)$ est au  plus 2.
Ceci ach\`eve la d\'emonstration du th\'eor\`eme. QED

\medskip

\noindent{\bf Remarque.} 
Du  th\'eor\`eme  \ref{te2} on d\'eduit imm\'ediatement, au moyen du  lemme \ref{trivial},
une majoration de l'ordre du groupe $\Coker(\alpha)$.

\section{Diff\'erentielles}\label{differentials}

\subsection{Une remarque g\'en\'erale sur les diff\'erentielles dans les suites spectrales}

Rappelons le cadre g\'en\'eral pour la suite spectrale des foncteurs compos\'es.
Soient  $\sA$, $\sB$, $\sC$ des cat\'egories ab\'eliennes. Supposons que
$\sA$ and $\sB$ ont suffisamment d'injectifs.
Soient $G:\sA\to\sB$ and $F:\sB\to\sC$  des foncteurs additifs exacts \`a gauche.
Supposons que $G$ envoie les objets injectifs sur des objets $F$-acycliques.
Alors pour tout objet $B\in{\rm Ob}(\sA)$ on a la suite spectrale
\begin{equation}
E_{2}^{pq}=(RF^p)(R^qG)B\Rightarrow R^{p+q}(FG)B.
\label{ss2}
\end{equation}
Soient 
$$\partial_{p,q} : (R^pF)(R^qG)B\ \lra \ (R^{p+2}F)(R^{q-1}G)B$$
les applications canoniques dans cette suite spectrale.

 Soit 
 \begin{equation}
0\to A\to B\to C\to 0 \label {ABC}
\end{equation}
une suite exacte dans $\sA$.
Par application des foncteurs d\'eriv\'es droits
de $G$ on obtient une longue suite exacte dans $\sB$.
En la tronquant on obtient pour tout $q\geq 1$
une suite exacte
\begin{equation}
0\to B_1\to (R^{q-1}G)C\to (R^q G)A\to B_2\to 0, \label{long}
\end{equation}une application surjective $s:(R^{q-1}G)B\to B_1$,
et une application injective $i:B_2\to (R^q G)B$. 
Soit $\partial: (R^pF)B_2\to (R^{p+2}F)B_1$ l'homorphisme
de connexion
d\'efini par (\ref{long}). Soit 
$$s_*=(R^{p+2}F)(s):(R^{p+2}F)(R^{q-1}G)B\to (R^{p+2}F)B_1$$
l'application induite par $s$, et, de fa{\c c}on analogue, soit
$$i_*=(R^pF)(i):(R^pF)B_2\to (R^pF)(R^q G)B$$
l'application induite par $i$.

\ble \label{te3}

On a  $\partial=s_*\partial_{p,q}i_*$.
\ele
{\it D\'emonstration.}  Soit
$$0\to A^\cdot\to B^\cdot\to C^\cdot\to 0$$ 
une suite exacte de r\'esolutions injectives de
 $A$, resp. $B$, resp. $C$.
Soient  $a_n:A^n\to A^{n+1}$ les diff\'erentielles
dans $A^\cdot$, et de m\^eme dans $B^\cdot$ et $C^\cdot$. 
On a le diagramme commutatif
$$\begin{array}{ccccccccc}
&& A^{q-1}/\Im(a_{q-2})&\to& B^{q-1}/\Im(b_{q-2})&\to& C^{q-1}/\Im(c_{q-2})&
\to &0\\
&&\downarrow&&\downarrow&&\downarrow&&\\
0&\to& \Ker(a_q)&\to& \Ker(b_q)&\to& \Ker(c_q)&& .
\end{array}$$
En appliquant le lemme du serpent on obtient la suite exacte
$$(R^{q-1}G)A\to (R^{q-1}G)B\to (R^{q-1}G)C\to
(R^{q}G)A\to (R^{q}G)B\to (R^{q}G)C.$$
En la tronquant on obtient
 (\ref{long}).  Par une chasse au diagramme on v\'erifie que la suite
(\ref{long}) est  \'equivalente \`a la  2-extension
\begin{equation}
0\to (R^{q-1}G)B\to B^{q-1}/b_{q-2}(B^{q-2}) \to \Ker(b_q) \to
(R^{q}G)B\to 0 \label{can}\end{equation}
tir\'ee en arri\`ere via
 $i:B_2\to (R^qG)B$ et pouss\'ee en avant via
$s:(R^{q-1}G)B\to B_1$. Par d\'efinition, l'application canonique
$\partial_{p,q}$ est l'homomorphisme de connexion 
$$(R^pF)(R^qG)B\lra(R^{p+2}F)(R^{q-1}G)B$$
d\'efini par  (\ref{can}), 
donc $s_*\partial_{p,q}i_*=\partial$. QED

\subsection{Applications au groupe de Brauer}

De la suite de Kummer  (\ref{kum}) on tire la 2-extension
de $\Ga$-modules
\begin{equation}
0\to\Pic(\ov X)/\Pic(\ov X)[n]\to\Pic(\ov X)\to\H^2_\et(\ov X,\mu_n)
\to\Br(\ov X)[n]\to 0, \label{2ext}
\end{equation}
o\`u la seconde fl\`eche est d\'efinie par la 
 multiplication
par $n$ sur  $\Pic(\ov X)$. 

\bpr \label{pr2}
Le diagramme suivant commute :
$$ \begin{array}{ccccccccc}
\Br(\ov X)[n]^{\Ga} & \buildrel{\partial}\over{\lra} &   
\H^2(k,\Pic(\ov X)/\Pic(\ov X)[n])\\
 \downarrow& & \uparrow \\
\Br(\ov X)^{\Ga}  & \buildrel{\beta}\over{\lra} &  \H^2(k,\Pic(\ov X))
\end{array}$$
Ici  $\partial$ est l'homomorphisme de connexion d\'efini par {\rm (\ref{2ext})},
et les fl\`eches verticales sont les applications naturelles  \'evidentes.
\epr
{\it D\'emonstration.} 
Dans le cadre de (\ref{ss2}) et du lemme \ref{te3}, soit
$\sA$ la cat\'egorie des faisceaux \'etales sur $X$, soit
$\sB$ la cat\'egorie des  $\Ga$-modules continus discrets, et soit
$\sC$ la cat\'egorie des groupes ab\'eliens.
Soit  $G=\pi_*$, o\`u $\pi:X\to\Spec(k)$
est le morphisme structural. Soit
$F(M)=M^\Ga$. Soit $A=\mu_{n,X}$, $B=C=\G_{m,X}$.
Pour  (\ref{ABC}),  prenons la suite de Kummer (\ref{kum}).
Prenons $p=0$ et $q=2$.
La suite exacte (\ref{long}) associ\'ee est pr\'ecis\'ement la suite
 (\ref{2ext}).
Il reste \`a appliquer le lemme \ref{te3} : 
la fl\`eche verticale de gauche dans le diagramme est  $i_{*}$,
la fl\`eche horizontale inf\'erieure est
$\beta = \partial_{0,2}$, 
et la fl\`eche verticale de droite est $s_{*}$. QED

\medskip

De la suite exacte
$$0\to \Pic^0(\ov X)\to\Pic(\ov X)\to\NS(\ov X)\to 0$$
on tire facilement la suite exacte
$$
0\to \Pic^0(\ov X)/\Pic^0(\ov X)[n]\to
\Pic(\ov X)/\Pic(\ov X)[n]\to\NS(\ov X)/\NS(\ov X)[n]\to 0. \label{av}
$$
Le sous-groupe  divisible   $\Pic^0(\ov X)\subset\Pic(\ov X)$ 
est contenu dans le noyau de
 $\Pic(\ov X)\to\H^2_\et(\ov X,\mu_n)$, donc de (\ref{2ext})
on tire la 2-extension de $\Ga$-modules
\begin{equation}
0\to\NS(\ov X)/\NS(\ov X)[n]\to\NS(\ov X)\to\H^2_\et(\ov X,\mu_n)
\to\Br(\ov X)[n]\to 0, \label{ns-ext}
\end{equation}
o\`u la seconde fl\`eche est induite par la multiplication par $n$
sur $\NS(\ov X)$.

\bco \label{c}
Le diagramme suivant commute :
$$ \begin{array}{ccccccccc}
\Br(\ov X)[n]^{\Ga} & \buildrel{\partial}\over{\lra} &
\H^2(k,\NS(\ov X)/\NS(\ov X)[n])\\
 \downarrow& & \uparrow \\
\Br(\ov X)^{\Ga}  & \buildrel{\beta}\over{\lra} &  \H^2(k,\Pic(\ov X))
 \end{array}$$
Dans ce diagramme,  $\partial$ est l'homomorphisme de connexion  d\'efini par {\rm (\ref{ns-ext})},
et les fl\`eches verticales sont les applications naturelles \'evidentes.
\eco
{\it D\'emonstration.} C'est une cons\'equence imm\'ediate de la proposition  \ref{pr2}. QED

\bco \label{c1}
Le diagramme suivant commute :
$$ \begin{array}{ccccccccc}
\Br^0(\ov X)^{\Ga} & \buildrel{\partial}\over{\lra} &
\H^2(k,N^1)\\
 \downarrow& & \uparrow \\
\Br(\ov X)^{\Ga}  & \buildrel{\beta}\over{\lra} &  \H^2(k,\Pic(\ov X))
 \end{array}$$
Dans ce diagramme,  $\partial$ est l'homomorphisme de connexion  d\'efini par {\rm (\ref{miau})},
et les fl\`eches verticales sont les applications naturelles \'evidentes.
\eco
{\it D\'emonstration.}
Soit $n$ un entier positif non nul divisible par  
l'exposant  $\nu_0$ de $\NS(\ov X)_{\rm tors}$.
Pour un tel  $n$ la suite exacte (\ref{ns-ext}) 
se lit
$$
0\to N^1 \to \NS(\ov X)\to\H^2_\et(\ov X,\mu_n)
\to\Br(\ov X)[n]\to 0,
$$
l'application  $N^1\to\NS(\ov X)$ \'etant induite par la
multiplication par $n$ sur $\NS(\ov X)$.
\'Ecrivons  $n=\prod_\ell n_{\ell}$, 
o\`u $n_{\ell}$ est une puissance du nombre premier $\ell$.

Soit  $P_{\ell}=\NS(\ov X)\{\ell\}$,
et soit $\Im(P_{\ell})$ l'image de  $P_\ell$ par l'application compos\'ee
$$\NS(\ov X) \to \H^2_{\et}(\ov X, \Z_{\ell}(1)) \to \H^2_\et(\ov X,\mu_{n_{\ell}}).$$
On a le diagramme commutatif suivant de
$\Ga$-modules,
dont les lignes sont exactes :
$$\begin{array}{ccccccccccc}
0&\to&N^1&\lra&N^1\otimes\Q&\to& 
\oplus_{\ell} (H^2_{\ell}\otimes_{\Z_\ell} \Q_\ell/\Z_\ell) &\to&
\Br^0(\ov X)&\to& 0\\
&&||&& \uparrow& &  \uparrow&&  \uparrow&  & \\
0&\to&N^1&\buildrel{\times n}\over{\lra}  &N^1 &\to& \oplus_{\ell} H^2_{\ell}/n_{\ell} &\to&
\Br^0(\ov X)[n]&\to& 0\\
 & & ||  &&   ||   & &  \downarrow&&  \downarrow&  & \\
0&\to&N^1&\buildrel{\times n}\over{\lra}  &N^1 &\to&  \oplus_{\ell}\H^2_\et(\ov X,\mu_{n_{\ell}})/\Im(P_{\ell}) &\to& \Br(\ov X)[n]&\to& 0\\
& & ||  &&   \uparrow  &&  \uparrow&& || &  & \\
0&\to&N^1&  \buildrel{\times n}\over{\lra}    & \NS(\ov X)&\to& \H^2_\et(\ov X,\mu_n) &\to& \Br(\ov X)[n]&\to& 0.
\end{array}
$$
La suite exacte de la premi\`ere ligne est obtenue en tensorisant
(\ref{tate}) avec $\Q_\ell/\Z_\ell$  puis en prenant la somme directe
sur tous les premiers $\ell$.
Pour construire la suite exacte de la seconde ligne, on tensorise
(\ref{tate}) avec  $\Z/n_\ell$, puis on prend la somme directe
sur tous les premiers $\ell$.
La fl\`eche {\it verticale} $N^1\to N^1\otimes\Q$ envoie
 $x$ sur  $x\otimes\frac{1}{n}$. Toutes les autres fl\`eches
 verticales sont les fl\`eches  naturelles \'evidentes.

En utilisant ce diagramme, on d\'eduit du corollaire \ref{c} que la restriction
de l'application compos\'ee 
$$\Br(\ov X)^\Ga \buildrel{\beta}\over{\lra}\H^2(k,\Pic(\ov X))\lra \H^2(k,N^{1})$$
\`a
$\Br^0(\ov X)^\Ga\subset\Br(\ov X)^\Ga$
est l'homomorphisme de connexion d\'efini par
 (\ref{miau}), la 2-extension sup\'erieure dans le  grand diagramme ci-dessus.
 QED

\section{D\'emonstration du th\'eor\`eme principal via les   cycles transcendants}\label{two}

Dans tout ce paragraphe, $X$ est une vari\'et\'e projective, lisse
et g\'eom\'etriquement int\`egre sur un corps $k$ de
caract\'eristique nulle.
 
\subsection{R\'eseaux de cycles alg\'ebriques et de cycles transcendants}

Soit $\ell$ un nombre premier. 
Pour 
 $M$ un  $\Z_\ell$-module, on note  $M^*=\Hom_{\Z_\ell}(M,\Z_\ell)$.
Dans le diagramme commutatif 
 d'accouplements $\Ga$-\'equivariants
$$\begin{array}{ccccc}
N^1_\ell& \times& N_{1,\ell} & \to& \Z_\ell\\
\downarrow&&\downarrow&&||\\
H^2_\ell&\times&H^{2d-2}_\ell&\to&\Z_\ell 
\end{array}$$
vu au paragraphe \ref{1.1}, 
 les fl\`eches verticales sont 
 {\it injectives} 
(th\'eor\`emes de Matsusaka et de Lieberman),
  En outre, l'accouplement inf\'erieur induit
 des isomorphismes
$H^2_\ell=(H^{2d-2}_\ell)^*$ et $H^{2d-2}_\ell=(H^2_\ell)^*$. 
En utilisant la suite exacte  
(\ref{x6}) et la remarque subs\'equente, on voit que la fl\`eche compos\'ee
$$N^1_\ell\to H^2_\ell\tilde\lra(H^{2d-2}_\ell)^*
\to N_{1,\ell}^*$$
est une application injective de conoyau
 $D\{\ell\}$,
o\`u $D$ est le groupe ab\'elien  fini d\'efini en (\ref{x4}).
En particulier, cette application est un isomorphisme si
$\ell$ ne divise pas l'ordre  $\delta$ de $D$.

Comme on a vu au paragraphe   \ref{1.2}, le sous-groupe   $N^1_\ell\subset H^2_\ell$
est primitif. Par contre  Koll\'ar (voir 
\cite{Voisin}, Thm. 14) a montr\'e que $N_{1,\ell}$ n'est pas forc\'ement un sous-groupe
primitif de
  $H^{2d-2}_\ell$. Voici comment rem\'edier \`a cet \'etat de choses.
Tout d'abord, si $\ell$ ne divise pas $\delta$, 
l'application naturelle $(H^{2d-2}_\ell)^*\to N_{1,\ell}^*$ est surjective, donc
$N_{1,\ell}$ est primitif dans $H^{2d-2}_\ell$. 
Pour tout
  $\ell$, on d\'efinit  le $\Ga$-module $M_\ell$ 
comme le satur\'e de  $N_{1,\ell}$ dans $H^{2d-2}_\ell$. En d'autres termes,
$$M_\ell\, =\ H^{2d-2}_\ell\ \cap (\ N_{1,\ell}\otimes_{\Z_\ell}\Q_\ell)
\ \subset \ H^{2d-2}_\ell\otimes_{\Z_\ell}\Q_\ell.$$
Si $\ell$ ne divise pas $\delta$, alors $M_\ell=N_{1,\ell}$.
Nous d\'efinissons ensuite le $\Ga$-module $M$ comme le sous-groupe de $N_1\otimes\Q$
form\'e des \'el\'ements qui,  pour tout premier $\ell$,
s'envoient dans $M_\ell$ par l'application naturelle 
$$N_1\otimes\Q\lra N_1\otimes\Q_\ell\cong
M_\ell\otimes_{\Z_\ell}\Q_\ell.$$
On a donc $N_{1} \subset M \subset N_1\otimes\Q$.
On obtient aussi une forme bilin\'eaire 
 $\Ga$-\'equivariante   $N^1\times M\to\Q$. 
Par tensorisation avec $\Z_\ell$ pour chaque premier  $\ell$ 
on voit que c'est en fait une forme bilin\'eaire enti\`ere
$$N^1\times M\,\to\,\Z$$
qui prolonge l'accouplement d'intersection sur
 $N^1\times N_1$.
Cela donne la suite exacte de $\Ga$-modules
\begin{equation}
0\to N^1\to \Hom(M,\Z)\to E\to 0, \label{m1}
\end{equation}
qui d\'efinit le $\Ga$-module fini $E$,
et pour chaque $\ell$ cela donne la suite exacte
\begin{equation}
0\to N^1_\ell\to M_\ell^*\to E\{\ell\}\to 0. \label{m2}
\end{equation}
Le $\Ga$-module 
  $E$ est un sous-$\Ga$-submodule of $D$,
et  $D/E=\Hom(M/N_{1},\Q/\Z)$. Donc  
$|D/E|=|M/N_1|$.
Notons que si $d=2$, c'est-\`a-dire si $X$ est une surface, 
alors $N_{1,\ell}=N^1_{\ell} \subset H^2_{\ell}$ est primitif, donc
  $M=N_{1}$ et  $D=E$.

Soit $S_\ell\subset H^2_\ell$ 
l'orthogonal de
$N_{1,\ell}$ (ou de  $M_\ell$)
par rapport au cup-produit.
Soit $T_\ell\subset H^{2d-2}_\ell$ 
l'orthogonal de
  $N^1_\ell$ par rapport au cup-produit.
  En dualisant la suite exacte de $\Gamma$-modules,
  libres  et de type fini comme  $\Z_{\ell}$-modules,
  $$0\to T_\ell\to H^{2d-2}_\ell\to (N^1_\ell)^*\to 0$$
on obtient la suite exacte
de $\Gamma$-modules,
   libres  et de type fini comme  $\Z_{\ell}$-modules,
\begin{equation}
0\to N^1_\ell\to H^2_\ell\to T_\ell^*\to 0.\label{m4}
\end{equation}
Ceci d\'efinit une identification canonique
 $T_\ell^*=B_\ell$,
o\`u $B_\ell$ est le module de Tate du groupe de Brauer d\'efini
au paragraphe \ref{1.2}.

Puisque  $M_\ell  \subset H^{2d-2}_\ell$  est un sous-groupe primitif,
le cup-produit donne 
la suite exacte suivante :
\begin{equation}
0\to S_\ell\to H^2_\ell\to M_\ell^*\to 0.\label{m3}
\end{equation}
L'application compos\'ee
  $N^1_\ell  \subset H^{2}_\ell \tilde\to (H^{2d-2}_\ell)^* \to  M_\ell^* \to (N_{1,\ell})^*$
est injective. Ainsi $ S_\ell \cap N^1_\ell=0$.
En utilisant  (\ref{m4}) and (\ref{m3}) on voit que pour tout premier $\ell$
 on a des isomorphismes canoniques de
 $\Ga$-modules
$$E\{\ell\}= M_\ell^*/N^1_\ell=H^2_\ell/(N^1_\ell\oplus S_\ell)=
T_\ell^*/S_\ell.$$
On a donc une suite exacte naturelle
\begin{equation}
0 \to M_\ell^*/N^1_\ell \to (S_\ell\otimes_{\Z_\ell}\Q_\ell/\Z_\ell) 
\oplus (N^1_\ell\otimes_{\Z_\ell}\Q_\ell/\Z_\ell)
\to  H^2_\ell\otimes_{\Z_\ell}\Q_\ell/\Z_\ell \to 0. \label{mm}
\end{equation}
Nous ferons usage du diagramme commutatif suivante de
 $\Ga$-modules, dont les colonnes et les lignes sont exactes : 
\begin{equation}\begin{array}{ccccccccc}
&&0&&0&&&&\\
&&\downarrow&&\downarrow&&&&\\
0&\to&M_\ell^*/N^1_\ell&\to&S_\ell\otimes_{\Z_\ell}\Q_\ell/\Z_\ell&
\to& \Br^0(\ov X)\{\ell\}&\to &0\\
&&\downarrow&&\downarrow&&||&&\\
0&\to& N^1_\ell\otimes_{\Z_\ell}\Q_\ell/\Z_\ell&\to& H^2_\ell\otimes_{\Z_\ell}\Q_\ell/\Z_\ell&
\to& \Br^0(\ov X)\{\ell\}&\to &0\\
&&\downarrow&&\downarrow&&&&\\
&&M_\ell^*\otimes_{\Z_\ell}\Q_\ell/\Z_\ell& =& M_\ell^*\otimes_{\Z_\ell}\Q_\ell/\Z_\ell&&&&\\
&&\downarrow&&\downarrow&&&&\\
&&0&&0&&&&
\end{array}\label{imp}
\end{equation}
La ligne m\'ediane, respectivement  la colonne m\'ediane,
est  la suite exacte
(\ref{m4}), 
respectivement  la suite exacte (\ref{m3}), tensoris\'ee avec  $\Q_\ell/\Z_\ell$.
Le reste du diagramme se d\'eduit de  (\ref{mm}).

En prenant la somme directe sur tous les premiers  $\ell$,  on d\'eduit de~(\ref{imp})
l'\'equivalence des 2-extensions
$$\begin{array}{ccccccccccc}
0&\to&N^1&\to&\Hom(M,\Z)&\to&\oplus(S_\ell\otimes_{\Z_\ell}\Q_\ell/\Z_\ell)&
\to& \Br^0(\ov X)&\to &0\\
&&||&&\downarrow&&\downarrow&&||&&\\
0&\to&N^1&\to& N^1\otimes\Q&\to& \oplus(H^2_\ell\otimes_{\Z_\ell}\Q_\ell/\Z_\ell)&
\to& \Br^0(\ov X)&\to &0
\end{array}
$$
L'extension sup\'erieure
est le produit de Yoneda des 
1-extensions de $\Ga$-modules
(\ref{m1}) et
\begin{equation}
0\to E\to\oplus(S_\ell\otimes_{\Z_\ell}\Q_\ell/\Z_\ell)
\to\Br^0(\ov X)\to 0.\label{m7}
\end{equation}
Notons $\partial_1:\Br^0(\ov X)^\Ga\to\H^1(k, E)$
et  $\partial_2:\H^1(k, E)\to \H^2(k,N^1)$ les diff\'erentielles
d\'efinies par ces 1-extensions.

\bpr \label{d2} 
L'application compos\'ee
$$\Br^0(\ov X)^{\Ga}\hookrightarrow 
\Br(\ov X)^{\Ga} \buildrel{\beta}\over{\lra} 
\H^2(k,\Pic(\ov X)) \to \H^2(k,N^1)$$
co\"{\i}ncide, au signe pr\`es, avec l'application compos\'ee
$$\Br^0(\ov X)^{\Ga}\buildrel{\partial_1}\over{\lra}\H^1(k,E)
\buildrel{\partial_2}\over{\lra} \H^2(k,N^1).$$
En particulier l'image de
 $\beta(\Br^0(\ov X)^{\Ga})$
dans $\H^2(k,N^1)$ est annul\'ee par l'exposant de  $E$.
\epr
{\it D\'emonstration.} 
Nous avons vu que 
(\ref{miau}) est \'equivalent au produit de Yoneda de
 (\ref{m1}) et  (\ref{m7}), 
 la proposition r\'esulte donc du corollaire
 \ref{c1}. QED

\subsection{Majorations, II}

Soit $\gamma$ l'ordre du groupe fini
$\oplus_{\ell} \H^3_{\et}(\ov X,\Z_{\ell}(1))_{\rm tors}$, et soit
  $\gamma_{0}$ son exposant. 
Soit $\varepsilon_{0} $ l'exposant du groupe fini $E$ 
d\'efini en (\ref{m1}).
L'entier $\varepsilon_{0}$ divise $\delta_0$, 
qui est l'exposant du groupe fini
  $D$  d\'efini en (\ref{x4}).  
Rappelons que $\nu$ est l'ordre du groupe fini  $\NS(\ov X)_{\rm tors}$,
et que  $\nu_0$ est son exposant.

Si $d=2$, c'est-\`a-dire si $X$ est une surface, alors
$\H^2_{\et}(\ov X,\Z_{\ell}(1))_{\rm tors}$ est dual de
$\H^{3}_{\et}(\ov X,\Z_{\ell}(1))_{\rm tors}$, donc $\gamma=\nu$
et $\gamma_{0}=\nu_{0}$. 
Dans ce cas $N^1=N_1$, et on a l'accouplement bilin\'eaire sym\'etrique
$$N^1 \times N^1\to \Z,$$
dont le noyau est trivial. L'entier $\delta=|D|$ est alors la valeur absolue
du d\'eterminant de cet accouplement. Toujours dans ce cas,
les groupes $D$ et $E$ co\"{\i}ncident,  donc
  $\delta=\varepsilon$ et $\delta_{0} =\varepsilon_{0}$.

Sur un corps quelconque (de caract\'eristique z\'ero), nous avons le r\'esultat suivant.

\bthe\label{t}
Soit $X$ une vari\'et\'e projective, lisse et g\'eom\'etriquement int\`egre sur un corps
$k$ de caract\'eristique nulle, telle que $\H^1(X,O_{X})=0$. 
Supposons que l'application canonique $\H^3_{\et}(k,\ov k^*)\to \H^3_\et(X,\G_{m})$ 
est injective (ce qui est le cas si $X$ poss\`ede un $k$-point). Alors :

{\rm (i)} L'exposant du conoyau de
$$ \alpha : \Br(X)\to \Br(\ov X)^\Gamma$$
divise  $\gamma_{0}\varepsilon_{0}\nu_0$, et
l'ordre de  $\Coker(\alpha)$ divise $\gamma (\varepsilon_{0}\nu_{0})^{b_{2}-\rho}$.

{\rm (ii)} Si $X$  est une surface, l'exposant de $\Coker(\alpha)$ divise $\delta_{0}\nu_{0}^2$, 
et l'ordre de
 $\Coker(\alpha)$ divise $\nu (\delta_{0}\nu_{0})^{b_{2}-\rho}$.
\ethe
{\it D\'emonstration.} Sous nos hypoth\`eses, $\Coker(\alpha)=\Im(\beta)$ 
d'apr\`es la proposition \ref{bas}, il suffit donc d'estimer la taille de
 $\beta(\Br(\ov X)^\Ga)$. 
De (\ref{brgeneral}) on d\'eduit la suite exacte
$$0 \to \Br^0(\ov X)^{\Ga} \to \Br(\ov X)^{\Ga} \to \oplus_{\ell}
\H^3_{\et}(\ov X,\Z_{\ell}(1))_{\rm tors}^{\Ga}.$$
Ceci implique que  $|\beta(\Br(\ov X)^\Ga)|$ divise
$\gamma|\beta({\Br^0(\ov X)}^\Ga)|$, et l'exposant de
 $\beta(\Br(\ov X)^\Ga)$ divise le produit de $\gamma_{0}$
par
 l'exposant de $\beta({\Br^0(\ov X)}^\Ga)$. 

D'apr\`es la proposition \ref{d2}, le groupe
$\varepsilon_{0}.\beta({\Br^0(\ov X)}^\Ga)$ est un sous-groupe de
$$\Ker[\H^2(k, \Pic(\ov X)) \to \H^2(k,N^1)].$$
On a la suite exacte courte
$$\H^2(k,\Pic^0(\ov X)) \to \H^2(k, \Pic(\ov X)) \to \H^2(k,\NS(\ov X))$$
et la suite exacte courte
$$ \H^2(k,\NS(\ov X)_{\rm tors})  \to \H^2(k,\NS(\ov X)) \to \H^2(k,N^1).$$
L'hypoth\`ese  $\H^1(X,O_{X})=0$ implique 
$\Pic^0(\ov X)=0$, donc l'exposant de
$\varepsilon_{0}.\beta({\Br^0(\ov X)}^\Ga)$ divise $\nu _{0}$.
Ainsi l'exposant de  $\beta({\Br^0(\ov X)}^\Ga)$ divise $\varepsilon_{0}.\nu_{0}$.
Le groupe  $\beta({\Br^0(\ov X)}^\Ga)$ est un sous-quotient de $(\Q/\Z)^{b_{2}-\rho}$.
Le lemme \ref{trivial} donne alors une borne pour l'ordre de
$\beta({\Br^0(\ov X)}^\Ga)$, qui donne la borne annonc\'ee pour
l'ordre de
   $\beta(\Br(\ov X)^\Ga)=\Coker(\alpha)$.
L'\'enonc\'e pour une surface r\'esulte des faits g\'en\'eraux rappel\'es au d\'ebut
de ce paragraphe.
QED

\medskip

Lorsque le corps de base  $k$ est un corps de nombres,
nous pouvons \'enoncer un r\'esultat sans la restriction  $\H^1(X,O_{X})=0$. 

\bthe\label{t1}
 Soit $X$une vari\'et\'e projective, lisse et g\'eom\'etriquement int\`egre sur un corps
 de nombres $k$. Alors :

{\rm (i)} L'exposant du conoyau de
$$ \alpha : \Br(X)\to \Br(\ov X)^\Gamma$$
divise  $2\gamma_{0}\varepsilon_{0} \nu_0$, et il divise
$\gamma_{0}\varepsilon_{0} \nu_0$
si $k$ est totalement imaginaire.
L'ordre de $\Coker(\alpha)$ divise  $\gamma (2\varepsilon_{0} \nu_{0})^{b_{2}-\rho}$,
et il divise  $\gamma (\varepsilon_{0}\nu_{0})^{b_{2}-\rho}$  
si $k$ est totalement imaginaire.

{\rm (ii)} Si $X$ est une surface, 
l'exposant de $\Coker(\alpha)$ divise $2\delta_{0}\nu_{0}^2$, et il divise
$\delta_{0}\nu_{0}^2$ si $k$ est totalement imaginaire; l'ordre de  $\Coker(\alpha)$
divise  $\nu(2\delta_{0}\nu_{0})^{b_{2}-\rho}$, et il divise  $\nu(\delta_{0}\nu_{0})^{b_{2}-\rho}$
si $k$ est totalement imaginaire.

\ethe
{\it D\'emonstration.} Pour un corps de nombres $k$, on a $\H^3_{\et}(k,\ov k^*)=0$.
Si l'on suit la d\'emonstration du th\'eor\`eme \ref{t}, le r\'esultat provient
du fait que le groupe
$\H^2(k,\Pic^0(\ov X))$ est un groupe fini d'exposant 2,
et que ce groupe est nul si $k$ est totalement imaginaire.
C'est un fait g\'en\'eral pour les vari\'et\'es ab\'eliennes sur un corps de nombres  \cite[Thm. 6.26 (c), p.~92]{Milne}.
L'\'enonc\'e pour les surfaces se d\'eduit de l'\'enonc\'e g\'en\'eral comme dans la pr\'ec\'edente d\'emonstration. QED

\medskip

\noindent{\bf Remarque} On a l'isomorphisme
$$\NS(\ov X)_{\rm tors}=\oplus_{\ell}\H^2_{\et}(\ov X,\Z_{\ell}(1))_{\rm tors}.$$
La dualit\'e de Poincar\'e implique que les groupes finis
$\H^2_{\et}(\ov X,\Z_{\ell}(1))_{\rm tors}$ et
$\H^{2d-1}_{\et}(\ov X,\Z_{\ell}(d-1))_{\rm tors}$ sont duaux l'un de l'autre.

\medskip

  La proposition suivante apporte un compl\'ement utile au th\'eor\`eme  \ref{t1}.
  Pour un corps de nombres $k$ on note $k_v$ la compl\'etion de $k$
  en une place non archim\'edienne $v$ et $k_v^{\rm nr}$ l'extension maximale non
  ramifi\'ee de $k_v$. Pour $S$ un ensemble fini de places finies 
  de $k$ et $E$ un module galoisien fini, on note $\H^1_S(k,E)$
  le sous-groupe de $\H^1(k,E)$ form\'e des \'el\'ements non ramifi\'es
  en dehors de $S$, c'est-\`a-dire l'intersection, pour
  tous les  $v\notin S$, des noyaux des 
  applications de restriction naturelles $\H^1(k,E)\to\H^1(k_v^{\rm nr},E)$.

\bpr \label{unr}
Soit $X$une vari\'et\'e projective, lisse et g\'eom\'e\-tri\-quement in\-t\`e\-gre sur un corps
 de nombres $k$, avec bonne r\'eduction  en dehors d'un ensemble fini $S$
 de places finies de $k$. Soit $E$ le $\Ga$-module  fini d\'efini
 en {\rm (\ref{m1})}.
Soit $T_E$ l'ensembles des places finies  de  $k$ divisant 
l'ordre de $E$.  Alors
$$\partial_1(\Br^0(\ov X)^{\Ga})\subset\H^1_{S\cup T_E}(k,E).$$
\epr
{\it D\'emonstration.} 
  Soit  $I_v=\Gal(\ov k_v/k_v)$ le groupe d'inertie.
Il est bien connu que si $v$ est une place de bonne r\'eduction,
et si la caract\'eristique r\'esiduelle est diff\'erente de  $\ell$, 
alors l'action naturelle de  $I_v$ sur $\H^2_{\et}(\ov X, \mu_{\ell^m})$, 
$m\geq 1$, et donc sur $\H^2_{\et}(\ov X, \Z_{\ell}(1))$, est triviale
(th\'eor\`eme de changement de base lisse, voir \cite[Cor. VI.4.2]{EC}).
Via la suite de Kummer, ceci implique que $I_v$ agit trivialement
sur $\Br(\ov X)\{\ell\}$. Mais $I_v$ agit aussi trivialement sur
$S_\ell\subset H^2_\ell$,
donc la diff\'erentielle
$$\partial:\Br^0(\ov X)\{\ell\}^{I_v}\lra \H^1(k_v^{\rm nr}, E\{\ell\})$$
d\'efinie par la suite exacte de  $I_v$-modules
$$0\to E\{\ell\}\to S_\ell\otimes_{\Z_\ell}\Q_\ell/\Z_\ell\to 
\Br^0(\ov X)\{\ell\}\to 0,$$
(ligne sup\'erieure de (\ref{imp})) est  nulle.
Donc, pour $v\notin S\cup T_E$, l'image de $\Br^0(\ov X)^{\Ga}$ dans  $\H^1(k,E)$
est dans le noyau de l'application de restriction \`a $\H^1(k_v^{\rm nr},E)$. QED

\section{Applications aux surfaces}\label{surfaces}

\bpr \label{K3}
Soit  $X$ une surface $K3$  sur un corps $k$ de caract\'eristique nulle. Supposons 
que l'application $\H^3_{\et}(k,\ov k^*)\to \H^3_\et(X,\G_{m})$ est injective (c'est le cas
si $k$ est un corps de nombres ou si $X$ poss\`ede un $k$-point).
Soit $\alpha$  l'application naturelle $\Br(X)\to \Br(\ov X)^\Gamma$.
L'exposant de  $\Coker(\alpha)$ divise $\delta_0$, et l'ordre de ce groupe
divise $\delta_0^{b_{2}-\rho}$.
\epr
{\it D\'emonstration.} Pour $X$ une surface  $K3$,  on a $\H^1(X,O_{X})=0$, et
le groupe de N\'eron--Severi
 $\NS(\ov X)$ est sans torsion, donc $\nu=1$. L'\'enonc\'e est un cas particulier du
th\'eor\`eme  \ref{t} (ii).
QED

\medskip

\noindent{\bf Exemples} 

1. Soit $X\subset \P^3_k$ une surface quartique diagonale sur
un corps de caract\'eristique z\'ero.
Il est bien connu que l'on a 
 $\delta=64$ et  $\delta_0=8$, voir \cite{PS}.
 Ceci implique d\'ej\`a que tout \'el\'ement d'ordre impair de 
 $\Br(\ov X)^\Ga$  vient de $\Br(X)$.
En outre, $b_{2}=22$ et $\rho=20$. 
En suivant la d\'emonstration du th\'eor\`eme  \ref{t}, on voit que
$\Coker(\alpha)$ est un sous-quotient de  $(\Q_{2}/\Z_{2})^2$.
Puisque son exposant divise $8$,  $\Coker(\alpha)$
est isomorphe \`a un sous-groupe de  $(\Z/8)^2$.

2.  Soit  $k$ un corps de caract\'eristique z\'ero.
Supposons que  $X\subset\P^g_k$
est une surface $K3$
``tr\`es g\'en\'erale'', c'est-\`a-dire que $\NS(\ov X)\cong\Z$ 
et est engendr\'e par la classe d'une section hyperplane~$H$.
On a  $\delta=(H.H)=2g-2$. 
La suite exacte (\ref{x4})  se lit alors
$$0 \to \Z \to \Z \to D \to 0,$$
o\`u  $E=D=\Z/(2g-2)$ avec $\Ga$-action triviale.
La d\'emonstration du th\'eor\`eme \ref{t} donne une injection
$$\Coker(\alpha)\hookrightarrow (\Z/(2g-2))^{21}.$$ 
Comme on a  $\H^1(k,\Z)=0$, l'application $$\partial_1 : 
\H^1(k,D) \to \H^2(k,\Z)=\H^2(k,\NS(\ov X))=\H^2(k,\Pic(\ov X))$$
est injective. Si $k$ est un corps de nombres,
$S$ est l'ensemble des places finies de mauvaise r\'eduction de  $X$, et
$T$ l'ensemble des places 
divisant $2g-2$, alors la proposition \ref{unr}
donne un homomorphisme injectif
$$\Coker(\alpha)\hookrightarrow\H^1_{S\cup T}(k,\Z/(2g-2)).$$

\medskip

Pour $X$ le produit de deux courbes on a un r\'esultat similaire au
th\'eor\`eme  \ref{t} (ii) bien qu'ici  $\H^1(X,O_X)\not=0$.

\bpr\label{p}
Soit $X=C_1\times C_2$ le produit de deux courbes projectives, lisses
et g\'eom\'etriquement int\`egres sur un corps $k$ de caract\'eristique nulle.
Soit $J_{1}$, resp. $J_{2}$, la jacobienne de 
  $C_{1}$, resp.  $C_{2}$.
Supposons que $X$ poss\`ede un $k$-point. Alors :

{\rm (i)} L'exposant du conoyau de
$\Br(X)\to \Br(\ov X)^{\Gamma}$ divise $\delta_0$. 

{\rm (ii)} Si  $\Hom_{\ov k}(\ov J_{1}, \ov  J_{2})=0$, alors
l'application
$\Br(X)\to \Br(\ov X)^{\Gamma}$ est surjective.
\epr
{\it D\'emonstration.}
Rappelons des  faits connus (\cite[\S 1; Cor. 6.2]{M3}).
L'application naturelle de  $\Ga$-modules
$$
(p_{1}^*, p_{2}^*) : \Pic(\ov C_1)\oplus \Pic(\ov C_2) \lra\Pic(\ov X)
$$
est une injection scind\'ee, dont une r\'etraction est donn\'ee par le
choix d'un
 $k$-point  $M=(M_{1},M_{2}) \in X(k)=C_{1}(k) \times C_{2}(k)$.
Ceci induit un isomorphisme de $\Ga$-modules
$$
\Pic^0(\ov C_1)\oplus \Pic^0(\ov C_2)\tilde\lra\Pic^0(\ov X).
$$
Le conoyau de  $(p_{1}^*,p_{2}^*)$ est le
  $\Ga$-module  $\Hom_{\ov k-{\rm grp}}(\ov J_{1},\ov J_{2})$, 
  qui est libre et de type fini comme  groupe ab\'elien.
Il y a une suite exacte induite  de $\Ga$-modules
 libres et de type fini comme  groupes ab\'eliens
$$0 \to \NS(\ov C_{1}) \oplus \NS(\ov C_{2}) \to \NS(\ov X) \to 
\Hom_{\ov k-{\rm grp}}(\ov J_{1},\ov J_{2}) \to 0,$$
c'est-\`a-dire
$$0\to\Z\oplus\Z\to\NS(\ov X)\to\Hom_{\ov k-{\rm grp}}(\ov J_{1},\ov J_{2}) \to 0,$$
un scindage \'etant donn\'e par le point $M$.
En utilisant le $k$-point $M\in X(k)$ on obtient donc un diagramme commutatif
\begin{equation}
\begin{array}{ccc}
\H^2(k,\Pic^0(\ov X))&\to &\H^2(k,\Pic(\ov X))\\
\downarrow{\simeq} &&\downarrow \\
\H^2(k,\Pic^0(\ov C_{1})\oplus \Pic^0(\ov C_{2}))&\hookrightarrow &
\H^2(k,\Pic(\ov C_{1})\oplus \Pic(\ov C_{2})).
\end{array}\label{di10}
\end{equation}
L'injection dans la ligne inf\'erieure vient de la nullit\'e de $\H^1(k,\Z)$.

Nous suivons la d\'emonstration du th\'eor\`eme \ref{te1}.
D'apr\`es la proposition  \ref{bas}, for $i=1,2$,  nous avons des diagrammes
commutatifs de suites exactes
$$\begin{array}{cccccc}
\Br(X) & \to & \Br(\ov X)^\Ga & \buildrel{\beta_{X}}\over{\lra} &\H^2(k,\Pic(\ov X))\\
\downarrow & &\downarrow&&\downarrow\\
\Br(C_{i}) & \to & \Br(\ov C_{i})^\Ga&
\buildrel{\beta_{C_{i}}}\over{\lra} &\H^2(k,\Pic(\ov C_{i})).
\end{array}$$
Par le th\'eor\`eme de Tsen, on a  $\Br(\ov C_{i})=0$. Ainsi
$$\beta(\Br(\ov X)^\Ga) \subset \Ker[\H^2(k,\Pic(\ov X)) \to  
\H^2(k,\Pic(\ov C_{1})\oplus \Pic(\ov C_{2}))  ].$$     
Puisque  $\NS(\ov X)$ est sans torsion, on voit
(remarque apr\`es le th\'eor\`eme \ref{t1})  que
$\H^3_{\et}(\ov X,\Z_{\ell} (1))$ est sans torsion pour tout $\ell$, 
ce qui implique $\Br^0(\ov X)=\Br(\ov X)$.
D'apr\`es la proposition \ref{d2}, l'image de $\Br(\ov X)^\Ga=\Br^0(\ov X)^\Ga$
dans le groupe $\H^2(k,N^1)=\H^2(k,\NS(\ov X))$ est annul\'ee par $\delta_0$.
Ainsi tout \'el\'ement dans
$\delta_0.\beta(\Br(\ov X)^\Ga) \subset \H^2(k,\Pic(\ov X))$ vient de
 $\H^2(k,\Pic^0(\ov X))$ et a une image nulle dans
$\H^2(k,\Pic(\ov C_{1})\oplus \Pic(\ov C_{2}))$.
De (\ref{di10}) on d\'eduit  $\delta_0.\beta(\Br(\ov X)^\Ga)=0$.
Ceci \'etablit  (i). 

Si $\Hom_{\ov k-{\rm grp}}(\ov J_{1}, \ov  J_{2})=0$, alors
$ \NS(\ov C_{1}) \oplus \NS(\ov C_{2}) \tilde\lra  \NS(\ov X) $ et
$\delta=1$.   
QED

\medskip

\noindent{\bf Exemples}

1. Lorsque $C_1=E$ et  $C_2=E'$ sont des courbes elliptiques non isog\`enes,
sur   $\ov k$,  on comparera l'\'enonc\'e ci-dessus avec
   \cite[Prop. 3.3]{SZ2}.
 
2. Si $C_1=C_2$ est une courbe elliptique $E$ sans multiplication complexe sur
 $\ov k$, alors $\NS(\ov X)$ est un groupe ab\'elien libre de rang 3
 engendr\'e par les classes de
  $E\times \{0\}$,  de $\{0\}\times E$
et de la diagonale $\Delta$. On a  $\delta =2$, $b_{2}=6$, $\rho=3$,
donc le conoyau de 
  $\alpha:\Br(X)\to \Br(\ov X)^\Gamma$ est isomorphe \`a
  un sous-groupe de
  $(\Z/2)^3$. Dans \cite[Prop. 4.3]{SZ2} on trouvera un exemple
avec $\Coker(\alpha) \neq 0$.

\section{Vari\'et\'es ouvertes}\label{open}

Nous offrons ici une r\'eponse partielle \`a une question soulev\'ee par
 T. Szamuely.

\bpr \label{h1ouvert} 
Soit $k$ un corps de type fini sur $\Q$.
Soit $U$ une $k$-vari\'et\'e quasi-projective
et lisse. Le groupe  $\H^1_{\et}(\overline{U},\Q/\Z)^\Ga$ est fini.
\epr
{\it D\'emonstration.} Ceci est une reformulation d'un cas particulier
d'un r\'esultat de 
Katz et  Lang \cite[Thm. 1, p.~295]{KL}. 
On peut aussi l'\'etablir par la m\'ethode plus g\'en\'erale suivante.
On peut supposer $U/k$ g\'eom\'etriquement connexe.
Par le th\'eor\`eme d'Hironaka, il existe une $k$-vari\'et\'e
projective, lisse et  g\'eom\'e\-tri\-quement int\`egre $X$
qui contient $U$ comme ouvert dense.
Soit $Z=X\setminus U$, et soit $F \subset Z$ le lieu singulier de $Z$.
Soit $X^{0} = X \setminus F$ et soit  $Z^0=Z \setminus  F$. Puisque  $X^0$ et
$F^0$ sont lisses, les suites de localisation pour la cohomologie \'etale 
\`a coefficients finis et le th\'eor\`eme de puret\'e  donnent
des suites exactes de
$\Gamma$-modules 
$$0 \to \H^1_{\et}(\ov X^0,\Q_{\ell}/\Z_{\ell}) \to 
\H^1_{\et}(\ov U,\Q_{\ell}/\Z_{\ell}) \to 
\H^0(\ov Z^0,\Q_{\ell}/\Z_{\ell}(-1)).$$ 
Puisque  $F$ est de codimension au moins  2 dans $X$, l'inclusion
$X^0\to X$ induit un isomorphisme
$$\H^1_{\et}(\ov X,\Q_{\ell}/\Z_{\ell})= \H^1_{\et}(\ov X^0,\Q_{\ell}/\Z_{\ell}).$$
On a donc une suite exacte
$$0 \to \H^1_{\et}(\ov X,\Q_{\ell}/\Z_{\ell})^\Gamma \to 
\H^1_{\et}(\ov U,\Q_{\ell}/\Z_{\ell})^\Gamma \to 
\H^0(\ov Z^0,\Q_{\ell}/\Z_{\ell}(-1))^\Gamma.$$ 
La  $k$-vari\'et\'e lisse  $Z^0$ se d\'ecompose comme union disjointe
$Z^0=\cup_{i} Z_{i}$ de $k$-vari\'et\'es lisses connexes.
Soit $k_{i}$ la fermeture int\'egrale de  $k$ dans le corps de fonctions $k(Z_{i})$.
Choisissons un $k$-plongement $k_{i} \subset \ov k$.
Soit $\Gamma_{i}=\Gal(\ov k/k_{i})$.
Le groupe  $\H^0(\ov Z^0,\Q_{\ell}/\Z_{\ell}(-1))^\Gamma$ est la somme directe
des groupes
$\oplus_{i} (\Q_{\ell}/\Z_{\ell}(-1))^{\Gamma_{i}}$.  Chaque groupe
$(\Q_{\ell}/\Z_{\ell}(-1))^{\Gamma_{i}}$ est fini; de plus, il est nul pour
presque tout~$\ell$. De fait, cet \'enonc\'e se ram\`ene 
imm\'ediatement \`a l'\'enonc\'e suivant : si $k$ est un corps de
nombres et  $\Gamma=\Gal(\ov k/k)$, alors  $(\Q_{\ell}/\Z_{\ell}(-1))^{\Gamma}$ 
est fini, et nul pour presque tout $\ell$. On est donc ramen\'e \`a v\'erifier
que $\H^1_{\et}(\ov X,\Q_{\ell}/\Z_{\ell})^\Gamma$  est fini et nul pour presque tout $\ell$.
Ceci r\'esulte du th\'eor\`eme de changement de base propre et des conjectures
de Weil (cf. \cite[Thm. 1.5]{CTR}).
QED

\bthe \label{h2ouvert}
Soit  $k$ un corps de type fini sur  $\Q$ et soit $U$ une $k$-vari\'et\'e
quasi-projective, lisse et g\'eom\'etriquement int\`egre sur $k$.
Alors :

{\rm (i)}  Le quotient $\Br(\ov U)^{\Gamma}/\Im(\Br(U))$ est un groupe fini.

{\rm (ii)} Si $U$ est une  surface, et si la conjecture de Tate  $\ell$-adique 
pour les diviseurs vaut pour une compactification lisse de  $U$, 
  $\Br(\ov U)\{\ell\}^{\Gamma}$ is finite.

{\rm (iii)} Si la conjecture de Tate  $\ell$-adique 
pour les diviseurs vaut pour une compactification lisse de  $U$, 
et si de plus le module galoisien $\H^2_{\et}(\ov X,\Q_{\ell}(1))$
est semi-simple,
alors   $\Br(\ov U)\{\ell\}^{\Gamma}$ est fini.
\ethe
{\it D\'emonstration.} Nous suivons la d\'emonstration de la proposition
\ref{h1ouvert} et nous utilisons les m\^emes notations.
Les suites exactes de localisation pour la cohomologie \'etale
\`a coefficients finis et le th\'eor\`eme de puret\'e donnent
naissance \`a la suite exacte de $\Gamma$-modules
$$0 \to \Br(\ov X^{0}) \to   \Br(\ov U) \to  \H^1_{\et}(\ov F^{0},\Q/\Z).$$
Comme la codimension de $F$ dans $X$ est au moins 2, le th\'eor\`eme de puret\'e
pour le groupe de Brauer montre que l'application de restriction
  $$\Br(\ov X) \to \Br(\ov X^{0})$$
est un  isomorphisme. On a donc la suite exacte de $\Gamma$-modules
$$0 \to \Br(\ov X) \to   \Br(\ov U) \to  \H^1_{\et}(\ov F^{0},\Q/\Z).$$
En prenant les invariants sous $\Gamma$, on obtient une suite exacte
$$0 \to \Br(\ov X)^{\Gamma} \to   \Br(\ov U)^{\Gamma} \to  
\H^1_{\et}(\ov F^{0},\Q/\Z)^{\Gamma}.$$
D'apr\`es la proposition \ref{h1ouvert}, le groupe  $\H^1_{\et}(\ov F^{0},\Q/\Z)^{\Gamma}$
est fini. Les \'enonc\'es (ii) and (iii) sont alors une cons\'equence de la finitude de 
 $ \Br(\ov X)\{\ell\}^{\Gamma}$,
laquelle vaut sous les hypoth\`eses de    (ii) ou (iii), cf.
\cite[Prop. 4.1]{CS}.

Par fonctorialit\'e on a un diagramme commutatif de suites exactes
$$\begin{array}{ccccccccccccccc}
0 & \to &  \Br(\ov X)^{\Gamma} & \to  &  \Br(\ov U)^{\Gamma} & \to  &  
{\H}^1_{\et}(\ov F^{0},\Q/\Z)^{\Gamma}  \\
&& \uparrow &&\uparrow&&   \\
&   & \Br(X) & \to  & \Br(U)  &  &  . \end{array}$$
D'apr\`es le th\'eor\`eme
  \ref{te1}, le quotient $\Br(\ov X)^{\Gamma} /\Im(\Br(X))$ est fini.
D'apr\`es la proposition \ref{h1ouvert}, le groupe  $\H^1_{\et}(\ov F^{0},\Q/\Z)^{\Gamma}$
est fini. Ceci implique que le quotient   $ \Br(\ov U)^{\Gamma}/\Im(\Br(U))$
est aussi fini, ce qui est l'\'enonc\'e  (i). QED

\bigskip

\noindent{\bf Remarque} Nous ne savons pas si (i) vaut sur tout corps
  $k$ de caract\'eristique nulle.

\vskip2cm

\noindent
CNRS, UMR 8628, Math\'ematiques, B\^atiment 425, Universit\'e Paris-Sud,
F-91405 Orsay, France
\medskip

\noindent jlct@math.u-psud.fr

\bigskip

\bigskip

\noindent Department of Mathematics, South Kensington Campus,
Imperial College London,
SW7 2BZ England, U.K.

\smallskip

\noindent Institute for the Information Transmission Problems,
Russian Academy of Sciences, 19 Bolshoi Karetnyi,
Moscow, 127994 Russia
\medskip

\noindent a.skorobogatov@imperial.ac.uk

\end{document}